\documentclass[12pt]{amsart}
\usepackage{amssymb}
\usepackage{amsmath}
\usepackage{amscd}
\usepackage[all]{xy}
\theoremstyle{plain}
\newtheorem{theorem}{Theorem}[section]
\newtheorem{corollary}[theorem]{Corollary}
\newtheorem{proposition}[theorem]{Proposition}
\newtheorem{lemma}[theorem]{Lemma}

{\theoremstyle{remark}

\newtheorem{remark}[theorem]{Remark}}
{\theoremstyle{definition}
\newtheorem{definition}[theorem]{Definition}

\newtheorem{example}[theorem]{Example}}

\newcommand{\rom}{\renewcommand{\labelenumi}{{\rm (\roman{enumi})}}%
\renewcommand{\itemsep}{0pt}}

\newcommand{\N}{\mathbb{N}}
\newcommand{\Z}{\mathbb{Z}}
\newcommand{\C}{\mathbb{C}}
\newcommand{\T}{\mathbb{T}}
\newcommand{\e}{\varepsilon}
\newcommand{\cK}{{\mathcal K}}
\newcommand{\cL}{{\mathcal L}}
\newcommand{\cM}{{\mathcal M}}

\newcommand{\F}{{\mathcal F}}
\newcommand{\Cr}{{\mathcal C}}
\newcommand{\D}{{\mathcal D}}
\newcommand{\tu}{\widetilde{u}}
\newcommand{\tw}{\widetilde{w}}
\newcommand{\cA}{\mathcal{A}}
\newcommand{\tcA}{\widetilde{\mathcal{A}}}
\newcommand{\cB}{\mathcal{B}}
\newcommand{\sX}{\mathsf{X}}

\newcommand{\fC}{{\mathfrak{C}}}
\newcommand{\act}{\curvearrowright}
\newcommand{\cO}{{\mathcal O}}
\newcommand{\cOs}{{\mathcal O}^{\rm{st}}}
\newcommand{\acO}{{\mathcal O}^{\rm{ alg}}}
\newcommand{\tcO}{{\widetilde{\mathcal O}}}
\newcommand{\cT}{{\mathcal T}}
\newcommand{\tcT}{{\widetilde{\mathcal T}}}
\newcommand{\cI}{{\mathcal I}}
\newcommand{\cJ}{{\mathcal J}}
\newcommand{\ip}[2]{\langle{#1},{#2}\rangle}

\newcommand{\s}[2]{s{\scriptstyle ({#1})}_{#2}}
\newcommand{\bS}[2]{S{({#1})}_{#2}}
\newcommand{\mat}[4]{\biggl(\!\begin{array}{cc}
{#1}&{#2}\\{#3}&{#4}\end{array}\!\biggr)}

\newcommand{\cD}{{\mathcal D}}
\newcommand{\tcD}{\widetilde{\mathcal D}}
\newcommand{\tphi}{\widetilde{\varphi}}
\newcommand{\tpi}{\widetilde{\pi}}

\newcommand{\Ca}{$C^*$-al\-ge\-bra }
\newcommand{\CA}{$C^*$-al\-ge\-bra}
\newcommand{\shom}{$*$-ho\-mo\-mor\-phism }
\newcommand{\shoms}{$*$-ho\-mo\-mor\-phisms }
\newcommand{\Cc}{$C^*$-cor\-re\-spon\-dence }

\newcommand{\Csa}{$C^*$-sub\-al\-ge\-bra }
\newcommand{\CsA}{$C^*$-sub\-al\-ge\-bra}

\DeclareMathOperator{\Aut}{Aut}

\DeclareMathOperator{\coker}{coker}

\DeclareMathOperator{\spa}{span}

\DeclareMathOperator{\cspa}{\overline{span}}

\begin{document}
\title[Construction of actions on Kirchberg algebras]
{A construction of actions on Kirchberg algebras 
which induce given actions on their K-groups.}
\author[Takeshi KATSURA]{Takeshi KATSURA}
\address{Department of Mathematics, 
Hokkaido University, Kita 10, Nishi 8, 
Kita-Ku, Sapporo, 060-0810, JAPAN}
\email{katsura@math.sci.hokudai.ac.jp}
\date{}

\keywords{$C^*$-algebra; Kirchberg algebra; group action; $K$-theory; 
Cuntz-Krieger algebra; group module.}

\subjclass[2000]{Primary 46L55; Secondary 46L40, 46L80}

\begin{abstract}
We prove that every action of 
a finite group all of whose Sylow subgroups are cyclic 
on the $K$-theory of a Kirchberg algebra 
can be lifted to an action on the Kirchberg algebra. 
The proof uses a construction of Kirchberg algebras 
generalizing the one of Cuntz-Krieger algebras, 
and a result on modules over finite groups. 
As a corollary, 
every automorphism of the $K$-theory of a Kirchberg algebra 
can be lifted to an automorphism of the Kirchberg algebra 
with same order. 
\end{abstract}

\maketitle

\setcounter{section}{-1}

\section{Introduction}

In this paper, 
we consider the problem of lifting a given action on the $K$-theory 
of a Kirchberg algebra $\cA$ to an action on $\cA$. 
We note that in this paper 
a Kirchberg algebra means 
a simple, separable, nuclear, purely infinite 
\Ca in the UCT class. 
First we state the problem more precisely. 

By the theorem of Kirchberg and Phillips in \cite{Ki,Ph}, 
Kirchberg algebras $\cA$ are classified 
by their $K$-theory $K_*(\cA)$ 
(see Subsection~\ref{ssec:KirAlg} for the definition). 
They also showed that 
the natural homomorphism $K_*\colon \Aut(\cA)\to \Aut(K_*(\cA))$ 
is surjective for every Kirchberg algebra $\cA$. 
An action of a group $\Gamma$ on a Kirchberg algebra $\cA$ 
induces an action of $\Gamma$ on $K_*(\cA)$ 
by composing the map $K_*$. 
The natural question asks whether 
we can get every action of $\Gamma$ on $K_*(\cA)$ in this way. 
This question can be restated as follows: 

\medskip
\noindent
{\bfseries Lifting Problem.} 
Let $\Gamma$ be a group 
and $\cA$ be a Kirchberg algebra. 
Does every action of $\Gamma$ on $K_*(\cA)$ 
lift to one on $\cA$? 
\medskip

To the author's best knowledge, 
no counterexample to this problem has been found so far. 
If one fixes a Kirchberg algebra $\cA$ 
and considers the lifting problem 
for arbitrary groups $\Gamma$, 
then one leads to the following splitting problem: 

\medskip
\noindent
{\bfseries Splitting Problem.} 
Does the surjection $K_*\colon \Aut(\cA)\to \Aut(K_*(\cA))$ split?
\medskip

This problem has been solved affirmatively 
for very few Kirchberg algebras $\cA$. 
We give results 
(Propositions~\ref{Prop:split1}, \ref{Prop:split2}, \ref{Prop:split3})
on the splitting problem 
as corollaries of our main theorem stated below. 
In this paper, 
we consider, instead of the splitting problem, 
the problem to find a group $\Gamma$ 
for which the lifting problem has an affirmative answer 
for arbitrary Kirchberg algebras. 
The group $\Gamma=\Z$ satisfies this property 
because the map $K_*\colon \Aut(\cA)\to \Aut(K_*(\cA))$ 
is surjective. 
For a finite group $\Gamma$,
as far as the author knows, 
the first result was due to 
Benson, Kumjian and Phillips 
who solved in \cite{BKP} the lifting problem affirmatively 
for $\Gamma=\Z/2\Z$ and for unital Kirchberg algebras $\cA$ 
in the Cuntz standard form. 
This result was extended by Spielberg 
who showed in \cite{Sp} that the lifting problem has an affirmative answer 
for $\Gamma=\Z/p\Z$ where $p$ is a prime number 
and for an arbitrary Kirchberg algebra. 
See also \cite{I} for another result on this problem 
using the Rohlin property. 
The following is the main theorem of this paper 
which extends the results in \cite{BKP} and \cite{Sp}. 

\medskip
\noindent
{\bfseries Theorem} (Theorem~\ref{Thm:Main}){\bfseries .} 
For a finite group $\Gamma$ all of whose Sylow subgroups are cyclic 
and an arbitrary Kirchberg algebra $\cA$, 
the lifting problem has an affirmative answer. 
\medskip

Since a finite cyclic group $\Z/n\Z$ satisfies the assumption 
in the main theorem, 
we get the following corollary. 

\medskip
\noindent
{\bfseries Corollary} (Corollary~\ref{Cor:Main}){\bfseries .} 
Let $\cA$ be a Kirchberg algebra. 
Then every automorphism of $K_*(\cA)$ 
can be lifted to an automorphism of $\cA$ 
with same order. 
\medskip

There are two ingredients of the proof of the main theorem. 
The first one is a construction of Kirchberg algebras $\cO_{A,B}$ 
using generators and relations 
coming from two matrices $A,B\in M_N(\Z)$. 
This construction can be considered as a generalization 
of the one of Cuntz-Krieger algebras introduced in \cite{CK}. 
We note that a similar construction can be found in \cite{D}. 
The good points of our construction are that 
we can construct an action of a certain group 
on our \Ca $\cO_{A,B}$ by permuting generators, 
and we can compute the $K$-theory $K_*(\cO_{A,B})$ 
as well as the induced action on it 
using two matrices $A,B$ (Proposition~\ref{Prop:KOAB}). 
This enables us to reduce the lifting problem to 
a problem on modules over finite groups (Theorem~\ref{Thm:constaction}). 
The other ingredient is the result in \cite{Ka5} 
which says that all modules have an affirmative answer to this problem 
for a finite group all of whose Sylow subgroups are cyclic 
(Theorem~\ref{Thm:permpre}). 
Combining these two results, 
we get the main theorem stated above. 

This paper is organized as follows. 
Section~\ref{Sec:Pre} is devoted to 
preparation of some notation and results we need 
in the rest of this paper. 
In Section~\ref{Sec:OAB}, 
we introduce the construction of a \Ca $\cO_{A,B}$ 
from two matrices $A,B\in M_N(\Z)$ 
satisfying a certain condition. 
We also define a group $\Gamma_{A,B}$ 
and its action on our \Ca $\cO_{A,B}$. 
We state three results on $\cO_{A,B}$ 
which say that we can compute 
the $K$-theory $K_*(\cO_{A,B})$ and 
the induced action of $\Gamma_{A,B}$ on it 
using two matrices $A,B\in M_N(\Z)$ 
(Proposition~\ref{Prop:KOAB}), 
$\cO_{A,B}$ is always separable, nuclear 
and in the UCT class (Proposition~\ref{Prop:OAB_sn}), and 
it is a Kirchberg algebra if $A,B$ satisfy certain conditions 
(Proposition~\ref{Prop:OAB_Ki}). 
Some parts of these results can be deduced from known results 
as explained in Section~\ref{Sec:Rel}. 
The author tries to make this paper self-contained 
as much as possible, 
and complete and direct proofs of all the three results 
are provided in the latter half of this paper. 
In Section~\ref{Sec:MainThm}, 
we prove the main theorem 
using the construction and results in Section~\ref{Sec:OAB} 
and the result on group modules from \cite{Ka5}. 
We also give several examples to which we can or cannot apply 
the main theorem, 
and show results on the splitting problem 
(Propositions~\ref{Prop:split1}, \ref{Prop:split2}, \ref{Prop:split3}). 
In Section~\ref{Sec:Rel}, 
we examine relations of our \CA s $\cO_{A,B}$ 
and the Cuntz-Krieger algebras $\cO_A$, 
the topological graph algebras $\cO(E_{A,B})$ 
and the Cuntz-Pimsner algebras $\cO_{\sX_{A,B}}$. 
We also explain which results in existing papers show 
which parts of the three results 
stated in Section~\ref{Sec:OAB}. 

In the rest of three sections, 
we give direct proofs of the three results on our \CA s $\cO_{A,B}$ 
in Section~\ref{Sec:OAB}. 
In Section~\ref{Sec:nuc}, 
we see that $\cO_{A,B}$ is always nuclear 
by examining the so-called core $\cO_{A,B}^{\T}$ of $\cO_{A,B}$. 
Using this analysis, 
we show in Section~\ref{Sec:pi} 
that $\cO_{A,B}$ becomes a Kirchberg algebra 
when $A,B\in M_N(\Z)$ satisfy certain three conditions. 
Finally in Section~\ref{Sec:K}, 
we compute the $K$-theory of $\cO_{A,B}$ 
and an action of the group $\Gamma_{A,B}$ on it 
induced by the action $\Gamma_{A,B}\act \cO_{A,B}$ 
defined in Section~\ref{Sec:OAB}. 
In the computation of the $K$-theories, 
we leave the proofs of two facts to two appendices. 
In Appendix~\ref{Sec:cA}, 
we prove Lemma~\ref{Lem:cA} 
which can be shown using \cite[Theorem~4.4]{Pi}. 
We try to make the proof very concrete, 
and to avoid $KK$-theory as much as possible. 
In Appendix~\ref{Sec:split}, 
we prove the last part of the proof of Proposition~\ref{Prop:exactseqs}. 
Although this part is not used in the proof of the main theorem 
(see Remark~\ref{Rem:Ker=0}), 
it is useful for examining several examples. 

\medskip

\noindent
{\bfseries Acknowledgments.} 
The author is grateful to the referee 
for careful reading. 
This work was partially supported by JSPS Research Fellow.

\section{Preliminaries}\label{Sec:Pre}

In this section, 
we prepare some notation and results. 

\subsection{Notation}
We denote by $\C$ the set of complex numbers, 
and by $\T\subset \C$ the group 
consisting of complex numbers with absolute values $1$. 
We denote by $\Z$ the abelian group of integers, 
and by $\N=\{0,1,\ldots \}\subset \Z$ 
the set of natural numbers. 
Let $\fC$ be the set of non-zero countable cardinalities 
which consists of positive integers $\{1,2,\ldots\}$
and the infinite countable cardinality $\infty$. 
For $N\in\fC$, 
$\{1,2,\ldots,N\}$ is a set with cardinality $N$. 
Here we promise that for $N=\infty$, 
$\{1,2,\ldots,N\}$ means 
the set of positive integers $\{1,2,\ldots\}$ 
similarly as 
$\{p_i\}_{i=1}^\infty$ usually means $\{p_1,p_2,\ldots\}$. 
We extend the addition on the set of positive integers 
to $\fC$ by setting $\infty + N=N+\infty =\infty$ 
for $N\in\fC$. 
Then the set $\{1,2,\ldots,N+N'\}$ has the same cardinality 
as $\{1,2,\ldots,N\}\amalg \{1,2,\ldots,N'\}$ 
for $N,N'\in\fC$. 

For a set $X$, 
we denote by $\ell^2(X)$ the Hilbert space 
whose complete orthonormal system is given 
by $\{\delta_{x}\}_{x\in X}$. 
For a Hilbert space $H$, 
we denote by $K(H)$ and $B(H)$ the \CA s of all compact operators 
and all bounded operators on $H$, respectively. 
For a subset $X$ of a \CA, 
$\spa X$ denotes the linear span of $X$, 
and $\cspa X$ denotes its closure. 

We denote by $\Aut(X)$ 
the group of automorphisms of a mathematical object $X$, 
such as a set, an abelian group or a $C^*$-al\-ge\-bra. 
An {\em action} $\Gamma\act X$ of a group $\Gamma$ 
on a mathematical object $X$ is 
a homomorphism from $\Gamma$ to $\Aut(X)$. 
For an action $\Gamma\act X$, 
the automorphism of $X$ defined by an element $\gamma\in \Gamma$ 
is often denoted by the same symbol $\gamma$. 
An abelian group $G$ with an action of $\Gamma$ 
is called a {\it $\Gamma$-module}. 
An isomorphism as $\Gamma$-modules means 
an isomorphism as abelian groups which is $\Gamma$-equivariant.

\subsection{Permutation presentations of modules}

Let us take $N\in\fC$. 
Let $\Z^N$ be the free abelian group 
whose basis is given by $\{e_i\}_{i=1}^N$. 
If an action $\Gamma\act \{1,2,\ldots,N\}$ is given, 
we can define an action $\Gamma\act \Z^N$ 
by $\gamma(e_i)=e_{\gamma(i)}$. 
Thus we get a $\Gamma$-module $\Z^N$. 
We call such a module a {\em countable permutation $\Gamma$-module}. 

\begin{definition}
A {\em countable permutation presentation} of a $\Gamma$-module $G$ is 
a $\Gamma$-equivariant exact sequence 
\[
\begin{CD}
0 @>>> F @>>> F @>>> G @>>> 0
\end{CD}
\]
where $F$ is a countable permutation $\Gamma$-module. 
\end{definition}

We are going to see that a countable permutation presentation is 
presented by a $N\times N$-matrix with integer entries for some $N\in\fC$. 

For a positive integer $N$, 
$M_N(\Z)$ denotes the algebra of 
all $N\times N$ matrices with integer entries. 
For $N=\infty$, 
$M_\infty(\Z)$ denotes the algebra of 
all infinite matrices $D=(D_{i,j})_{i,j=1}^\infty$ 
with integer entries satisfying that the set 
\[
\big\{j\in\{1,2,\ldots\}\ \big|\ D_{i,j}\neq 0\big\}
\]
is finite for all $i$. 

Let $N\in\fC$. 
A homomorphism $\varphi\colon \Z^N\to \Z^N$ is represented by 
an element $D\in M_N(\Z)$ 
so that $\varphi(e_i)=\sum_{j=1}^N D_{i,j}e_j$, 
and by this correspondence we can and will identify the set $M_N(\Z)$ 
and the set of all endomorphisms of $\Z^N$. 
Note that 
we consider $\Z^N$ as the set of row vectors via 
\[
(n_1,n_2,\ldots,n_N)\mapsto \sum_{i=1}^N n_ie_i \in \Z^N, 
\]
and think that $M_N(\Z)$ acts $\Z^N$ from right. 
We denote by $\N^N$ the subset of $\Z^N$ 
consisting of elements 
in the form $f=\sum_{i=1}^N n_ie_i\in \Z^N$ 
with $n_i\in\N$ for all $i$, 
and by $M_N(\N)$ the subset of $M_N(\Z)$ 
consisting of the matrices whose entries are in $\N$. 
The identity matrix of $M_N(\Z)$ 
is denoted by $I\in M_N(\Z)$. 

Let us take an action $\Gamma\act \{1,2,\ldots,N\}$. 
This action makes $\Z^N$ a countable permutation $\Gamma$-module. 
We set 
\[
M_N(\Z)^\Gamma
:=\{D\in M_N(\Z)\mid 
\text{$D_{i,j}=D_{\gamma(i),\gamma(j)}$ for all $i,j$ and $\gamma\in\Gamma$}\}.
\]
It is easy to see that a matrix $D$ is in $M_N(\Z)^\Gamma$ 
if and only if the homomorphism $\Z^N\to \Z^N$ 
determined by $D$ is $\Gamma$-equivariant. 
For $D\in M_N(\Z)^\Gamma$, 
the action $\Gamma\act \Z^N$ induces 
actions of $\Gamma$ on the abelian groups $\ker D$ and $\coker D$. 
Thus $\ker D$ and $\coker D$ become $\Gamma$-modules. 
We will denote by $[\cdot]$ the natural surjection $\Z^N\to \coker D$. 
The following observation is easy to see. 

\begin{lemma}\label{Lem:permpreD}
Giving a countable permutation presentation of a $\Gamma$-module $G$ is 
same as giving 
$N\in\fC$, an action $\Gamma \act \{1,2,\ldots,N\}$, 
a matrix $D\in M_N(\Z)^\Gamma$ with $\ker D=0$ 
and an isomorphism $\coker D\cong G$ as $\Gamma$-modules. 
\end{lemma}

\begin{example}\label{Ex1}
Let $\Gamma_1\cong \Z/2\Z$ be a cyclic group generated 
by $\sigma$ with $\sigma^2=1$. 
Let $G=\Z/3\Z$ be a $\Gamma_1$-module 
where an action of $\Gamma_1$ on $G$ is defined by $\sigma(g)=-g$. 
Then $G$ has a countable permutation presentation 
which is represented by $N=2$, an action of $\Gamma_1$ on $\{1,2\}$ 
defined by $\sigma(1)=2$ and $\sigma(2)=1$, 
a $2\times 2$-matrix 
\[
D=\mat{2}{-1}{-1}{2}\in M_2(\Z)^{\Gamma_1} 
\]
with $\ker D=0$, 
and an isomorphism $\coker D\cong G$ defined 
by $[e_i]\mapsto i$ for $i=1,2$. 
\end{example}

\begin{example}\label{Ex2}
Let $\Gamma_2\cong (\Z/2\Z)\times (\Z/2\Z)$ be a group 
generated by two elements $\sigma,\tau\in\Gamma_2$ 
with relations $\sigma^2=\tau^2=1$ and $\sigma\tau=\tau\sigma$. 
Let $G=(\Z/3\Z)^3$ be a $\Gamma_2$-module 
where an action of $\Gamma_2$ on $G$ is defined by 
\[
\sigma\big((a_1,a_2,a_3)\big)=(a_2, a_1, -a_1-a_2-a_3),\quad 
\tau\big((a_1,a_2,a_3)\big)=(a_3,-a_1-a_2-a_3, a_1), 
\]
for $a_1,a_2,a_3\in \Z/3\Z$. 
This $\Gamma_2$-module $G$ has a countable permutation presentation 
which is represented by $N=4$, 
an action of $\Gamma_2$ on $\{1,2,3,4\}$ defined by 
\begin{align*}
\sigma &\colon 1\mapsto 2,\quad 2\mapsto 1,\quad 3\mapsto 4,\quad 4\mapsto 3,\\
\tau &\colon 1\mapsto 3,\quad 2\mapsto 4,\quad 3\mapsto 1,\quad 4\mapsto 2,
\end{align*}
a matrix 
\[
D=\left(\begin{array}{cccc}
2&-1&-1&-1\\
-1&2&-1&-1\\
-1&-1&2&-1\\
-1&-1&-1&2
\end{array}\right)\in M_4(\Z)^{\Gamma_2}, 
\]
with $\ker D=0$ and an isomorphism $\coker D\cong G$ defined by 
\[
\Big[\sum_{i=1}^4n_ie_i\Big]\mapsto 
(n_1+n_2+n_3,n_1+n_2+n_4,n_1+n_3+n_4). 
\]
\end{example}

\begin{example}\label{Ex3}
The same formula as in Example~\ref{Ex2} defines 
an action of $\Gamma_2$ on $G'=(\Z/4\Z)^3$. 
We can show that this $\Gamma_2$-module $G'$ 
has no countable permutation presentations 
(see \cite[Example~2.6]{Ka5}). 
\end{example}

In \cite{Ka5}, 
we prove the following. 

\begin{theorem}[{\cite[Proposition~1.3 and Theorem~1.4]{Ka5}}]
\label{Thm:permpre}
Let $\Gamma$ be a finite group all of whose Sylow subgroups are cyclic. 
Then every countable $\Gamma$-module has 
a countable permutation presentation. 
\end{theorem}

Note that if a finite group $\Gamma$ has a non-cyclic Sylow subgroup, 
then the countable $\Gamma$-module called the augmentation ideal of $\Gamma$ 
has no countable permutation presentations 
(see \cite[Proposition~2.8]{Ka5}).

\subsection{\boldmath{$K$}-groups and partial unitaries}\label{ssec:partunit}

Let $\cA$ be a \CA . 
For definitions and results of $K$-groups $K_0(\cA)$ and $K_1(\cA)$, 
we consult the book \cite{Bl}. 
We denote by $[p]$ the element in $K_0(\cA)$ 
defined by a projection $p\in \cA$. 
For a projection $p\in \cA$, 
we say that $u$ is a {\em partial unitary with} $u^0=p$ 
if $u^*u=uu^*=p$. 
We denote by $\tcA$ the minimal unital \Ca containing $\cA$. 
For a partial unitary $u\in\cA$ with $u^0=p$, 
we denote by $\tu\in\tcA$ 
the unitary $u+(1-p)$. 
The element in $K_1(\cA)$ defined by 
this unitary $\tu$ 
is denoted by $[u]$. 
For a partial unitary $u\in\cA$ with $u^0=p$, 
we define $u^n\in\cA$ for $n\in\Z$ by 
\[
u^n=
\begin{cases}
u^n & \text{for $n>0$,}\\
p & \text{for $n=0$,}\\
(u^*)^{-n} & \text{for $n<0$.}\\
\end{cases}
\]
Then $u^n\in\cA$ is also a partial unitary with $(u^n)^0=p$ 
and satisfies $\widetilde{u^n}=\tu^n$ for $n\in\Z$. 
For a finite family $\{u_i\}_i$ of mutually orthogonal partial unitaries 
in $\cA$, 
the element $u=\sum_iu_i\in \cA$ is a partial unitary. 
The family $\{\tu_i\}_i$ of unitaries in $\tcA$ is mutually commutative, 
and we have $\tu=\prod_i\tu_i$. 
Hence we get $[u]=\sum_i[u_i]$ in $K_1(\cA)$.

\subsection{Kirchberg algebras}\label{ssec:KirAlg}

\begin{definition}
A {\em Kirchberg algebra} is a simple, separable, nuclear, purely infinite 
\Ca in the UCT class. 
\end{definition}

For a detailed definition and results of Kirchberg algebras, 
we consult the book \cite{RS}. 
We remark that 
in \cite[Definition~4.3.1]{RS} or some literatures, 
Kirchberg algebras are not assume to be in the UCT class. 

For a Kirchberg algebra $\cA$, 
we define $K_*(\cA)$ by 
\[
K_*(\cA)=\begin{cases}
\big(K_0(\cA),K_1(\cA)\big)&\text{if $\cA$ is stable,}\\
\big(K_0(\cA),[1_\cA],K_1(\cA)\big)&\text{if $\cA$ has a unit $1_\cA$.}
\end{cases}
\]
Recall that a Kirchberg algebra is stable 
if and only if it is non-unital (\cite{Z}). 

Let $G_0,G_1,G'_0$ and $G'_1$ be 
countable abelian groups. 
An isomorphism of the two pairs $(G_0,G_1)$ 
and $(G_0',G_1')$ is defined to be 
a pair $\varphi_*=(\varphi_0,\varphi_1)$ 
of isomorphisms $\varphi_i\colon G_i\to G'_i$ for $i=0,1$. 
If such an isomorphism exists, 
we write $(G_0,G_1)\cong (G_0',G_1')$. 
Let us take $g\in G_0$ and $g'\in G'_0$. 
An isomorphism of $(G_0,g,G_1)$ 
and $(G_0',g',G_1')$ is defined to be 
a pair $\varphi_*=(\varphi_0,\varphi_1)$ 
of isomorphisms $\varphi_i\colon G_i\to G'_i$ for $i=0,1$ 
with $\varphi_0(g)=g'$. 
If such an isomorphism exists, 
we write $(G_0,g,G_1)\cong (G_0',g',G_1')$. 
With these preparations, 
we can state the celebrated classification theorem due to 
Elliott, Kirchberg, Phillips and R\o rdam 
as follows 
(see \cite[Proposition~4.3.3, Theorem~8.4.1]{RS}). 

\begin{theorem}\label{Thm:KP}
Let $(G_0,G_1)$ be a pair of countable abelian groups. 
Then there exists a unique stable Kirchberg algebra $\cA$ 
with $K_*(\cA)\cong (G_0,G_1)$, 
and for each $g\in G_0$ 
there exists a unique unital Kirchberg algebra $\cA$ 
with $K_*(\cA)\cong (G_0,g,G_1)$. 
\end{theorem}

On the line of the proof of the theorem above, 
the natural homomorphism 
$\Aut(\cA)\to \Aut(K_*(\cA))$ is shown to be surjective 
for a Kirchberg algebra $\cA$. 
By composing this surjection, 
an action of a group $\Gamma$ on a Kirchberg algebra $\cA$ 
induces an action of $\Gamma$ on the $K$-theory $K_*(\cA)$ of $\cA$. 
In particular, 
the two groups $K_0(\cA)$ and $K_1(\cA)$ become $\Gamma$-modules.

\section{The \Ca $\cO_{A,B}$}\label{Sec:OAB}

In this section, 
we construct a \Ca $\cO_{A,B}$, a group $\Gamma_{A,B}$ 
and an action $\Gamma_{A,B}\act \cO_{A,B}$ 
from two matrices $A\in M_N(\N)$ and $B\in M_N(\Z)$. 
We also state three results on $\cO_{A,B}$ 
which will be proven in the latter half of this paper. 

\begin{definition}
Let $N\in\fC$. 
For $A\in M_N(\N)$, 
we define a set $\Omega_{A}$ by 
\[
\Omega_{A}
:=\big\{(i,j)\in \{1,2,\ldots,N\}\times \{1,2,\ldots,N\}\ \big|\ 
A_{i,j}\geq 1\big\}. 
\]
For each $i\in\{1,2,\ldots,N\}$, 
we define a set $\Omega_{A}(i)\subset \{1,2,\ldots,N\}$ by 
\[
\Omega_{A}(i)
:=\big\{j\in \{1,2,\ldots,N\}\ \big|\ (i,j)\in \Omega_{A}\big\}. 
\]
\end{definition}

Note that by definition $\Omega_{A}(i)$ is finite for all $i$. 

\begin{definition}\label{Def:OAB}
Let us take $N\in\fC$, 
$A\in M_N(\N)$ and $B\in M_N(\Z)$. 
We define a \Ca $\cO_{A,B}$ to be the universal \Ca 
generated by mutually orthogonal projections $\{p_i\}_{i=1}^N$, 
partial unitaries $\{u_i\}_{i=1}^N$ with $u_i^0=p_i$, 
and partial isometries $\{\s{n}{i,j}\}_{(i,j)\in \Omega_{A},n\in\Z}$ 
satisfying the relations 
\begin{enumerate}
\rom
\item $\s{n}{i,j}u_{j}=\s{n+A_{i,j}}{i,j}$ and 
$u_{i}\s{n}{i,j}=\s{n+B_{i,j}}{i,j}$ 
for all $(i,j)\in \Omega_{A}$ and $n\in\Z$, 
\item $\s{n}{i,j}^*\s{n}{i,j}=p_j$ 
for all $(i,j)\in \Omega_{A}$ and $n\in\Z$, 
\item 
$p_i=\sum_{j\in \Omega_{A}(i)}
\sum_{n=1}^{A_{i,j}}\s{n}{i,j}\s{n}{i,j}^*$ 
for all $i$. 
\end{enumerate}
\end{definition}

If there exists $i$ with $\Omega_{A}(i)=\emptyset$, 
then the condition (iii) says that $p_i=0$. 
Hence $\cO_{A,B}$ is isomorphic to $\cO_{A',B'}$ 
where $A',B'\in M_{N-1}(\Z)$ are obtained 
by eliminating the $i$-th rows and the $i$-th columns from $A,B$. 
By repeating this argument as many as possible, 
either we get $\cO_{A,B}=0$, or 
we can find $A'\in M_{N'}(\N)$ and $B'\in M_{N'}(\Z)$ 
which are isomorphic to corners of $A$ and $B$ 
such that $\Omega_{A'}(i)\neq\emptyset$ for all $i$ 
and $\cO_{A,B}\cong \cO_{A',B'}$ naturally. 
Hence without loss of generality, 
we may assume that $\Omega_{A}(i)\neq\emptyset$ for all $i$. 
In Definition~\ref{Def:OAB}, 
we use $B_{i,j}\in\Z$ only for $(i,j)\in \Omega_{A}$. 
Hence without loss of generality, 
we may assume that $(i,j)\not\in \Omega_{A}$ implies $B_{i,j}=0$. 
We summarize these assumptions 
into the following condition on $A,B\in M_N(\Z)$ 
for further reference; 
\begin{itemize}
\item[(0)] 
$A\in M_N(\N)$, $\Omega_{A}(i)\neq \emptyset$ for all $i$, 
and $B_{i,j}=0$ for $(i,j)\not\in \Omega_{A}$. 
\end{itemize}

\begin{remark}
As we did when defining graph algebras (see \cite{Ra}), 
we can change ``for all $i$'' in the condition (iii) 
to ``for all $i$ with $\Omega_{A}(i)\neq\emptyset$'' 
in order to get a meaningful \Ca $\cO_{A,B}$ 
for $A\in M_N(\N)$ and $B\in M_N(\Z)$ 
with $\Omega_{A}(i)=\emptyset$ for some $i$. 
One can also drop the assumption that 
$\Omega_{A}(i)$ is finite 
by changing the condition (iii) suitably 
(cf.\ \cite[Proposition~B.2]{Ka4}). 
For our purpose in this paper 
we do not need such a generality, 
and hence 
we only consider $A,B\in M_N(\Z)$ satisfying the condition (0). 
\end{remark}

\begin{definition}
Let us take $N\in\fC$, 
and $A,B\in M_N(\Z)$. 
We define a group $\Gamma_{A,B}$ by 
\[
\Gamma_{A,B}
:=\big\{\gamma\in \Aut(\{1,2,\ldots,N\})\ \big|\ 
\text{$A_{i,j}=A_{\gamma(i),\gamma(j)}$ and $B_{i,j}=B_{\gamma(i),\gamma(j)}$ 
for all $i,j$}\big\}. 
\]
\end{definition}

By definition, 
$\Gamma_{A,B}$ acts on $\{1,2,\ldots,N\}$ 
and we have $A,B\in M_N(\Z)^{\Gamma_{A,B}}$. 
If we have $A,B\in M_N(\Z)^{\Gamma}$ 
for an action of some group $\Gamma$ on $\{1,2,\ldots,N\}$, 
then there exists a unique homomorphism $\Gamma\to \Gamma_{A,B}$ 
such that the action of $\Gamma$ is the composition of this homomorphism 
and the action of $\Gamma_{A,B}$. 

\begin{definition}
Let $A,B\in M_N(\Z)$ satisfy the condition (0) above. 
We define an action $\Gamma_{A,B}\act \cO_{A,B}$ by 
\[
\gamma(p_i)=p_{\gamma(i)}, \quad
\gamma(u_i)=u_{\gamma(i)}\quad \text{and} \quad
\gamma(\s{n}{i,j})=\s{n}{\gamma(i),\gamma(j)} 
\]
for the generators $\{p_i,u_i,\s{n}{i,j}\}$ of the \Ca $\cO_{A,B}$ 
and $\gamma\in \Gamma_{A,B}$. 
\end{definition}

It is routine to check that 
the definition above is well-defined. 
This action induces actions of $\Gamma_{A,B}$ on 
$K_0(\cO_{A,B})$ and $K_1(\cO_{A,B})$. 
On the other hand, 
since $I-A,I-B\in M_N(\Z)^{\Gamma_{A,B}}$, 
we get actions of $\Gamma_{A,B}$ 
on $\ker (I-A)$, $\coker (I-A)$, $\ker (I-B)$ and $\coker (I-B)$. 
We obtain the following 
whose proof can be found in Section~\ref{Sec:K}. 

\begin{proposition}\label{Prop:KOAB}
There exist $\Gamma_{A,B}$-equivariant isomorphisms 
\begin{align*}
K_0(\cO_{A,B})
&\cong \coker (I-A) \oplus \ker (I-B),\\
K_1(\cO_{A,B})
&\cong \coker (I-B)\oplus \ker (I-A) 
\end{align*}
under which 
$[p_i]\in K_0(\cO_{A,B})$ and $[u_i]\in K_1(\cO_{A,B})$ 
correspond to $[e_i]\in \coker (I-A)$ and 
$[e_i]\in \coker (I-B)$ respectively for every $i$. 
\end{proposition}

\begin{proof}
This follows from Proposition~\ref{Prop:exactseqs} and its proof. 
\end{proof}

\begin{corollary}\label{Cor:KOAB}
Let $A,B\in M_N(\Z)$ satisfy the condition (0) above. 
Let $\Gamma\act\{1,2,\ldots,N\}$ be an action 
with $A,B\in M_N(\Z)^\Gamma$ which naturally induces 
the actions of $\Gamma$ 
on $\ker (I-A)$, $\coker (I-A)$, $\ker (I-B)$ and $\coker (I-B)$. 
Then we have the natural action $\Gamma\act\cO_{A,B}$ 
such that the isomorphisms in Proposition \ref{Prop:KOAB} 
are $\Gamma$-equivariant. 
\end{corollary}

\begin{remark}
For $A,B\in M_N(\Z)$ satisfying the condition (0), 
the \Ca $\cO_{A,B}$ is unital if and only if $N<\infty$, 
and in this case the unit is $\sum_{i=1}^N p_i\in\cO_{A,B}$ 
(see Subsection~\ref{ssec:CK}). 
Thus for $N<\infty$, 
$[1_{\cO_{A,B}}]\in K_0(\cO_{A,B})$ 
corresponds to $\big[\sum_{i=1}^N e_i\big]\in \coker (I-A)$ 
under the isomorphism in Proposition \ref{Prop:KOAB}. 
\end{remark}

We also get the following two propositions 
which will be proven 
in Sections~\ref{Sec:nuc}, \ref{Sec:pi} and \ref{Sec:K}. 

\begin{proposition}\label{Prop:OAB_sn}
For $A,B\in M_N(\Z)$ satisfying the condition {\rm (0)}, 
the \Ca $\cO_{A,B}$ is separable, nuclear 
and in the UCT class. 
\end{proposition}

\begin{proof}
The \Ca $\cO_{A,B}$ is separable 
because its generator is countable. 
It is nuclear by Proposition~\ref{Prop:OAB_nuc}, 
and in the UCT class by Proposition~\ref{Prop:OAB_UCT}. 
\end{proof}

\begin{proposition}\label{Prop:OAB_Ki}
If $A,B\in M_N(\Z)$ satisfy the condition {\rm (0)}
and the conditions 
\begin{enumerate}
\item $A\in M_N(\N)$ is irreducible, i.e. 
for every $i,j\in \{1,2,\ldots,N\}$ 
there exists a positive integer $n$ with $(A^n)_{i,j}\geq 1$, 
\item $A_{i,i}\geq 2$ and $B_{i,i}=1$ for every $i\in \{1,2,\ldots,N\}$, 
\end{enumerate}
then the \Ca $\cO_{A,B}$ 
is simple and purely infinite, 
and hence a Kirchberg algebra by Proposition~\ref{Prop:OAB_sn}. 
\end{proposition}

\begin{proof}
See Section~\ref{Sec:pi}. 
\end{proof}

\section{The proof of the main theorem}\label{Sec:MainThm}

In this section, 
we prove the main theorem (Theorem~\ref{Thm:Main}) 
using the results in the previous section and Theorem~\ref{Thm:permpre}. 
The following lemma is 
an equivariant version of 
\cite[Lemma~6.4]{Ka4}. 

\begin{lemma}\label{Lem:matrix}
Let $\Gamma$ be a finite group. 
For $N'\in\fC$, 
an action $\Gamma\act \{1,2,\ldots,N'\}$ 
and $A',B'\in M_{N'}(\Z)^\Gamma$, 
there exist $N\in\fC$, 
an action $\Gamma\act \{1,2,\ldots,N\}$ 
and $A,B\in M_{N}(\Z)^\Gamma$ 
satisfying the conditions {\rm (0)}, {\rm (1)}, {\rm (2)} 
in Section~\ref{Sec:OAB} 
such that there exist $\Gamma$-equivariant isomorphisms 
\begin{align*}
\ker(I-A)&\cong\ker A', &
\coker(I-A)&\cong\coker A', \\
\ker(I-B)&\cong\ker B', &
\coker(I-B)&\cong\coker B'. 
\end{align*}
\end{lemma}

\begin{proof}
In this proof, 
we denote by $I'\in M_{N'}(\N)$ the identity matrix of $M_{N'}(\N)$. 
We define $|A'|,|B'|\in M_{N'}(\N)$ 
by $|A'|_{i,j}=|A'_{i,j}|$ and $|B'|_{i,j}=|B'_{i,j}|$ 
for $i,j\in\{1,2,\ldots,N'\}$. 
We define $X\in M_{N'}(\N)$ by 
\[
X_{i,j}=\begin{cases}
1 & \text{if }|i-j|= 1,\\
0 & \text{if }|i-j|\neq 1,
\end{cases}
\]
and $X^\gamma\in M_{N'}(\N)$ by $(X^\gamma)_{i,j}=X_{\gamma(i),\gamma(j)}$ 
for $\gamma\in\Gamma$. 
We set $Y\in M_{N'}(\N)$ by 
\[
Y=|A'|+|B'|+I'+\sum_{\gamma\in\Gamma}X^{\gamma}. 
\]
Then we have $Y\in M_{N'}(\Z)^\Gamma$.
We define $A,B\in M_{2}(M_{N'}(\Z))$ by 
\[
A=\mat{2I'}{A'+Y}{I'}{I'+Y}, \quad
B=\mat{I'}{B'}{I'}{I'}.
\]
We set $N=N'+N'\in \fC$. 
Choose a bijection 
\[
\{1,2,\ldots,N'\}\amalg \{1,2,\ldots,N'\}\cong \{1,2,\ldots,N\} 
\]
and fix it. 
Using this bijection, 
we identify $\Z^{N'}\oplus \Z^{N'}\cong \Z^{N}$ and
$M_2(M_{N'}(\Z))\cong M_{N}(\Z)$, 
and regard $A,B\in M_{N}(\Z)$. 
Through this bijection, 
the action $\Gamma\act \{1,2,\ldots,N'\}$ 
induces an action $\Gamma\act \{1,2,\ldots,N\}$. 
We see that $A,B\in M_{N}(\Z)^\Gamma$ 
because all the entries of $A,B\in M_{2}(M_{N'}(\Z))$ 
are in $M_{N'}(\Z)^\Gamma$. 
By noticing $A'+|A'|\in M_{N'}(\N)$, 
one can easily check that $A,B\in M_{N}(\Z)$ 
satisfy the conditions (0) and (2). 
Since $X^{1}=X\in M_{N'}(\N)$ is irreducible, 
so is 
\[
\mat{I'}{I'}{I'}{X^{1}}\in M_{N}(\N). 
\] 
Hence $A\in M_{N}(\N)$ is irreducible. 
Thus $A,B\in M_{N}(\Z)$ also satisfy the condition (1). 

Note that the identity matrix $I\in M_{N}(\Z)$ corresponds to 
\[
\mat{I'}{0}{0}{I'}\in M_2(M_{N'}(\Z))
\]
via the identification $M_{N}(\Z)\cong M_2(M_{N'}(\Z))$. 
We have the equality 
\begin{align*}
I-A&=\mat{-I'}{-A'-Y}{-I'}{-Y}\\
&=\mat{I'}{I'}{0}{I'}\mat{A'}{0}{0}{-I'}\mat{0}{-I'}{I'}{Y}. 
\end{align*}
Since the left and the right matrices in the multiplication above 
define $\Gamma$-equivariant isomorphisms 
from $\Z^{N}$ to $\Z^{N}$, 
there exist $\Gamma$-equivariant isomorphisms 
\begin{align*}
&\ker(I-A)\cong 
\ker\mat{A'}{0}{0}{-I'}
\cong \ker A',\\
&\coker(I-A)\cong 
\coker\mat{A'}{0}{0}{-I'}
\cong \coker A'.
\end{align*}
Similarly, 
we have $\Gamma$-equivariant isomorphisms 
\[
\ker(I-B)\cong\ker B', \quad 
\coker(I-B)\cong\coker B'. 
\]
We are done. 
\end{proof}

\begin{proposition}\label{Prop:M->AB}
Let $\Gamma$ be a finite group, and 
$G_0,G_1$ be $\Gamma$-modules 
which have countable permutation presentations. 
Then there exist $N\in\fC$, 
an action $\Gamma\act \{1,2,\ldots,N\}$ 
and $A,B\in M_{N}(\Z)^\Gamma$ 
satisfying the conditions {\rm (0)}, {\rm (1)}, {\rm (2)} 
in Section~\ref{Sec:OAB} 
such that $\ker (I-A)=\ker (I-B)=0$, 
$\coker (I-A)\cong G_0$ and $\coker (I-B)\cong G_1$ as $\Gamma$-modules.

Moreover, for each $g\in G_0$ 
fixed by the action of $\Gamma$, 
we can find $f\in \N^N\subset \Z^N$ fixed by the action of $\Gamma$ 
such that $[f]\in \coker (I-A)$ corresponds to $g\in G_0$ 
under the isomorphism above. 
\end{proposition}

\begin{proof}
For $i=0,1$, 
Lemma~\ref{Lem:permpreD}
gives $N_i\in\fC$, 
an action $\Gamma\act\{1,2,\ldots, N_i\}$, 
a matrix $D_i\in M_{N_i}(\Z)^\Gamma$ with $\ker D_i=0$ 
and an isomorphism $\coker D_i\cong G_i$ as $\Gamma$-modules. 
Let $N'=N_0+N_1$ and choose a bijection 
\[
\{1,2,\ldots, N_0\}\amalg\{1,2,\ldots, N_1\}\cong 
\{1,2,\ldots, N'\}. 
\]
Through this bijection, 
two actions of $\Gamma$ on $\{1,2,\ldots, N_i\}$ for $i=0,1$ 
define an action of $\Gamma$ on $\{1,2,\ldots, N'\}$. 
Let $A',B'\in M_{N'}(\Z)^\Gamma$ be the images 
of $(D_0,I),(I,D_1)\in M_{N_0}(\Z)^\Gamma\oplus M_{N_1}(\Z)^\Gamma$ 
under the natural inclusion 
$M_{N_0}(\Z)^\Gamma\oplus M_{N_1}(\Z)^\Gamma\to M_{N'}(\Z)^\Gamma$ 
defined from the bijection. 
Then we have $\ker A'=\ker B'=0$, 
$\coker A'\cong G_0$ and $\coker B'\cong G_1$ as $\Gamma$-modules.
By Lemma~\ref{Lem:matrix}, 
we get $N\in\fC$, 
an action $\Gamma\act \{1,2,\ldots,N\}$ 
and $A,B\in M_{N}(\Z)^\Gamma$ 
satisfying the conditions (0), (1), (2) in Section~\ref{Sec:OAB} 
such that $\ker (I-A)=\ker (I-B)=0$, 
$\coker (I-A)\cong G_0$ and $\coker (I-B)\cong G_1$ as $\Gamma$-modules.

The isomorphism $\coker (I-A)\cong G_0$ gives us 
a surjection $\pi\colon \Z^N\to G_0$ 
such that $[f]\in \coker (I-A)$ corresponds to 
$\pi(f)\in G_0$ for each $f\in \Z^N$. 
The kernel of the surjection $\pi$ 
is isomorphic to $\Z^N$ 
because $\ker(I-A)=0$: 
\[
\begin{CD}
0 @>>> \Z^N @>I-A>> \Z^N @>\pi>> G_0 @>>> 0 
\end{CD}
\]
Take an element $g\in G_0$ fixed by the action of $\Gamma$. 
We can show that there exists 
$f\in \Z^N$ which is fixed by the action of $\Gamma$ 
and satisfies $\pi(f)=g$ 
because the obstruction of such a lifting is encoded 
in the group cohomology $H^1(\Gamma,\ker\pi)\cong H^1(\Gamma,\Z^N)$ 
which can be shown to vanish easily 
(see \cite[Remark~2.3 and Lemma~3.4~(1)]{Ka5}). 
Let us write $f=\sum_{i=1}^N n_ie_i$ for $n_i\in \Z$. 
If $n_i\geq 0$ for all $i$, 
this $f\in \N^N\subset \Z^N$ satisfies the desired condition. 
Suppose that there exists $i_0$ with $n_{i_0}<0$. 
Let $f_{i_0}\in \Z^N$ be the image of $e_{i_0}\in \Z^N$ under the map $A-I$. 
Then we have $f_{i_0}-e_{i_0}\in\N^N$ 
because $A\in M_N(\N)$ and $A_{i_0,i_0}\geq 2$. 
Set $f':=f-n_{i_0}\sum_{\gamma\in\Gamma}\gamma(f_{i_0})$. 
Then $f'$ is fixed by the action of $\Gamma$ 
and satisfies $\pi(f')=\pi(f)=g$ 
because $\pi(f_{i_0})=0$. 
If we write $f'=\sum_{i=1}^N n'_ie_i$ for $n'_i\in \Z$
then we have $n'_i\geq n_i$ for all $i$ 
and $n'_{i_0}\geq 0$. 
Repeating this argument, 
we get an element $f\in \N^N\subset \Z^N$ 
which is fixed by the action of $\Gamma$ 
and satisfies $\pi(f)=g$. 
We are done. 
\end{proof}

\begin{theorem}\label{Thm:constaction}
Let $\Gamma$ be a finite group, and $\cA$ be a Kirchberg algebra. 
An action $\Gamma\act K_*(\cA)$ 
lifts to an action $\Gamma\act \cA$ 
if the induced two $\Gamma$-modules $K_0(\cA)$ and $K_1(\cA)$ have 
countable permutation presentations. 
\end{theorem}

\begin{proof}
With the help of Theorem~\ref{Thm:KP}, 
it suffices to find 
\begin{itemize}
\item a stable Kirchberg algebra $\cA_s$ with an action of $\Gamma$ 
such that $K_*(\cA_s)\cong (G_0,G_1)$ $\Gamma$-equivariantly, and 
\item a unital Kirchberg algebra $\cA_u$ 
with an action of $\Gamma$ 
such that $K_*(\cA_u)\cong (G_0,g,G_1)$ $\Gamma$-equivariantly  
\end{itemize}
for $\Gamma$-modules $G_0,G_1$ having 
countable permutation presentations 
and $g\in G_0$ fixed by the action of $\Gamma$. 

Take $\Gamma$-modules $G_0,G_1$ having 
countable permutation presentations. 
Let $N\in\fC$, 
an action $\Gamma\act \{1,2,\ldots,N\}$ 
and $A,B\in M_{N}(\Z)^\Gamma$ be as 
in the conclusion of Proposition~\ref{Prop:M->AB}. 
Set $\cA:=\cO_{A,B}$ 
which is a Kirchberg algebra by Proposition~\ref{Prop:OAB_Ki}. 
Since $A,B\in M_{N}(\Z)^\Gamma$ 
we have an action $\Gamma\act \cA$, 
and by Corollary~\ref{Cor:KOAB} 
we get 
$K_i(\cA)\cong G_i$ as $\Gamma$-modules for $i=0,1$. 
Let $\cK$ be the \Ca of all compact operators 
on the separable infinite dimensional Hilbert space $\ell^2(\N)$. 
We define $\cA_s:=\cA\otimes \cK$ 
which is a stable Kirchberg algebra. 
The action $\Gamma\act \cA$ 
extends to an action $\Gamma\act \cA_s$ 
by acting $\cK$ trivially, 
and we have a $\Gamma$-equivariant isomorphism 
$K_*(\cA_s)\cong (G_0,G_1)$. 

Now take an element $g\in G_0$ fixed by the action of $\Gamma$. 
By Proposition~\ref{Prop:M->AB}, 
there exists $f\in \N^N\subset \Z^N$ fixed by the action of $\Gamma$ 
such that $[f]\in \coker (I-A)$ corresponds to $g\in G_0$. 
Let us denote $f=\sum_{i=1}^N n_ie_i$ 
with $n_i\in\N$. 
For each $n\in\N$, 
choose a projection $q_n\in \cK$ whose rank is $n$. 
Let us define 
$p=\sum_{i=1}^N p_i\otimes q_{n_i}\in \cA\otimes \cK=\cA_s$. 
Since $f\in\Z^N$ is fixed by the action of $\Gamma$, 
the projection $p\in \cA_s$ is also fixed by the action of $\Gamma$. 
Hence the action $\Gamma\act \cA_s$ globally fixes 
the unital Kirchberg algebra $\cA_u:=p\cA_s p$. 
By Proposition~\ref{Prop:KOAB}, 
the element $[p]\in K_0(\cA_s)$ 
corresponds to $[f]\in \coker (I-A)$ 
and hence to $g\in G_0$ 
by the isomorphisms $K_0(\cA_s)\cong \coker (I-A)\cong G_0$. 
Hence we have 
$K_*(\cA_u)\cong (G_0,g,G_1)$ $\Gamma$-equivariantly. 
We are done. 
\end{proof}

\begin{remark}
It seems to be possible to prove the theorem above 
using the construction in \cite{Sp}.
\end{remark}

\begin{theorem}\label{Thm:Main}
Let $\Gamma$ be a finite group all of whose Sylow subgroups are cyclic 
and $\cA$ be a Kirchberg algebra. 
Then every action $\Gamma \act K_*(\cA)$ 
lifts to an action $\Gamma \act \cA$. 
\end{theorem}

\begin{proof}
Combine Theorem~\ref{Thm:permpre} and Theorem~\ref{Thm:constaction}. 
\end{proof}

\begin{corollary}\label{Cor:Main}
Let $\cA$ be a Kirchberg algebra. 
Then every automorphism of $K_*(\cA)$ 
can be lifted to an automorphism of $\cA$ 
with same order. 
\end{corollary}

\begin{proof}
For an automorphism of $K_*(\cA)$ with infinite order, 
this follows from \cite{Ki,Ph}. 
For an automorphism of $K_*(\cA)$ with finite order $n$, 
this follows from Theorem~\ref{Thm:Main} 
because a finite cyclic group $\Z/n\Z$ satisfies 
the assumption of Theorem~\ref{Thm:Main}. 
\end{proof}

We are going to see some examples to which 
we can or cannot apply 
Theorem~\ref{Thm:constaction} or Theorem~\ref{Thm:Main}. 
A unital Kirchberg algebra $\cA$ is said to be 
{\em in the Cuntz standard form} 
if $[1_\cA]=0$ in $K_0(\cA)$. 
For $n=2,3,\ldots,\infty$, 
we denote by $\cOs_n$ the unital Kirchberg algebra 
in the Cuntz standard form with isomorphic $K$-groups 
as the Cuntz algebra $\cO_n$. 
Thus we have $K_*(\cOs_{n+1})=(\Z/n\Z,0,0)$ for $n<\infty$ 
and $K_*(\cOs_{\infty})=(\Z,0,0)$. 
Note that $\cOs_{n+1}\cong M_{n}(\cO_{n+1})$ for $n<\infty$.

\begin{example}
Let a finite group $\Gamma_1\cong \Z/2\Z$ and 
a $\Gamma_1$-module $G=\Z/3\Z$ be as in Example~\ref{Ex1}. 
In this example, 
we construct a unital Kirchberg algebras $\cA$ 
and an action $\Gamma_1\act \cA$ 
such that $K_0(\cA)\cong G$ as $\Gamma_1$-modules and $K_1(\cA)=0$ 
using the countable permutation presentation in Example~\ref{Ex1}. 
Note that $0\in G$ is the only element 
which is fixed by the action of $\Gamma_1$. 

Let us define an action $\Gamma_1\act\{1,2,3,4\}$ 
by $\sigma(1)=2$, $\sigma(2)=1$, $\sigma(3)=4$ and $\sigma(4)=3$.
Let $A,B\in M_4(\Z)^{\Gamma_1}$ be 
\[
A=
\left(\begin{array}{cccc}
2&0&3&0\\
0&2&0&3\\
1&0&2&1\\
0&1&1&2
\end{array}\right),\qquad 
B=
\left(\begin{array}{cccc}
1&0&1&0\\
0&1&0&1\\
1&0&1&0\\
0&1&0&1
\end{array}\right). 
\]
These two matrices are obtained by applying the argument 
in the proof of Lemma~\ref{Lem:matrix} 
for $A',B',Y\in M_2(\Z)^{\Gamma_1}$ defined by 
\[
A'=\mat{2}{-1}{-1}{2},\quad 
B'=\mat{1}{0}{0}{1},\quad
Y=\mat{1}{1}{1}{1}. 
\]
Although this $Y$ does not coincide with 
the one obtained from $A',B'$ 
as in the proof of Lemma~\ref{Lem:matrix}, 
we can still show that $A,B$ 
satisfies the conditions {\rm (0)}, {\rm (1)}, {\rm (2)}, 
and we have $\ker (I-A)=\ker (I-B)=\coker(I-B)=0$ 
and $\coker(I-A)\cong G$ as $\Gamma_1$-modules. 
Set $\cA:=\cO_{A,B}$ which is a unital Kirchberg algebra. 
Since $A,B\in M_4(\Z)^{\Gamma_1}$, 
we get an action $\Gamma_1\act \cA$ and 
we have $K_0(\cA)\cong G$ as $\Gamma_1$-modules and $K_1(\cA)=0$. 
By Theorem~\ref{Thm:KP}, 
$\cA$ is isomorphic to $\cOs_4$. 
\end{example}

\begin{example}
Similarly as in the previous example, 
we see that 
the \Ca $\cA=\cO_{A,B}$ 
is a unital Kirchberg algebra 
with an action $\Gamma_1\act \cA$ 
such that $K_0(\cA)\cong K_1(\cA)\cong G$ as $\Gamma_1$-modules 
where an action $\Gamma_1\act\{1,2,3,4\}$ 
is same as in the previous example and 
$A,B\in M_4(\Z)^{\Gamma_1}$ is defined by 
\[
A=
\left(\begin{array}{cccc}
2&0&3&1\\
0&2&1&3\\
1&0&2&2\\
0&1&2&2
\end{array}\right),\qquad 
B=
\left(\begin{array}{cccc}
1&0&2&-1\\
0&1&-1&2\\
1&0&1&0\\
0&1&0&1
\end{array}\right). 
\]
By Theorem~\ref{Thm:KP}, 
$\cA$ is isomorphic to $\cOs_4\otimes \cOs_4$. 
\end{example}

\begin{example}
Let a finite group $\Gamma_2\cong (\Z/2\Z)\times (\Z/2\Z)$ and 
a $\Gamma_2$-module $G=(\Z/3\Z)^3$ 
be as in Example~\ref{Ex2}. 
We construct a unital Kirchberg algebras $\cA$ 
and an action $\Gamma_2\act \cA$ 
such that $K_0(\cA)\cong G$ 
as $\Gamma_2$-modules and $K_1(\cA)=0$. 
Note that $0\in G$ is the only element 
fixed by the action of $\Gamma_2$. 

We define an action $\Gamma_2\act\{1,2,\ldots,8\}$ by 
\begin{align*}
\sigma &\colon 
1\mapsto 2,\quad 2\mapsto 1,\quad 3\mapsto 4,\quad 4\mapsto 3,\quad 
5\mapsto 6,\quad 6\mapsto 5,\quad 7\mapsto 8,\quad 8\mapsto 7,\\
\tau &\colon 
1\mapsto 3,\quad 2\mapsto 4,\quad 3\mapsto 1,\quad 4\mapsto 2, \quad 
5\mapsto 7,\quad 6\mapsto 8,\quad 7\mapsto 5,\quad 8\mapsto 6.
\end{align*}
We set $A,B\in M_8(\Z)^{\Gamma_2}$ by 
\[
A=
\left(\begin{array}{cccccccc}
2&0&0&0&3&0&0&0\\
0&2&0&0&0&3&0&0\\
0&0&2&0&0&0&3&0\\
0&0&0&2&0&0&0&3\\
1&0&0&0&2&1&1&1\\
0&1&0&0&1&2&1&1\\
0&0&1&0&1&1&2&1\\
0&0&0&1&1&1&1&2
\end{array}\right),\qquad 
B=
\left(\begin{array}{cccccccc}
1&0&0&0&1&0&0&0\\
0&1&0&0&0&1&0&0\\
0&0&1&0&0&0&1&0\\
0&0&0&1&0&0&0&1\\
1&0&0&0&1&0&0&0\\
0&1&0&0&0&1&0&0\\
0&0&1&0&0&0&1&0\\
0&0&0&1&0&0&0&1
\end{array}\right). 
\]
Then $\cA:=\cO_{A,B}$ is a unital Kirchberg algebra, 
and the action $\Gamma_2\act\{1,2,\ldots,8\}$ 
defines an action $\Gamma_2\act \cA$ 
so that $K_0(\cA)\cong G$ 
as $\Gamma_2$-modules and $K_1(\cA)=0$. 

We have $\Gamma_{A,B}\cong{\mathfrak S}_4$ 
where ${\mathfrak S}_4:=\Aut(\{1,2,3,4\})$. 
Thus we get an action ${\mathfrak S}_4\act \cA$ 
such that the action $\Gamma_2\act \cA$ is its restriction. 
The induced action ${\mathfrak S}_4\act K_0(\cA)\cong (\Z/3\Z)^3$ 
coincides with the natural action of ${\mathfrak S}_4$ on 
\[
\big\{(a_1,a_2,a_3,a_4)\in (\Z/3\Z)^4 \ \big|\ a_1+a_2+a_3+a_4=0\big\}
\cong (\Z/3\Z)^3. 
\]
Note that the group $\Gamma_2$ as well as ${\mathfrak S}_4$ 
does not satisfy the assumption of Theorem~\ref{Thm:Main}. 
Note also that the ${\mathfrak S}_4$-module $(\Z/3\Z)^3$ 
has a countable permutation presentation 
as we can see from the matrix $A$ above or 
the matrix $D$ in Example~\ref{Ex2}. 
\end{example}

\begin{example}
The author does not know how to construct 
a Kirchberg algebra $\cA$ 
and an action $\Gamma_2\act \cA$ 
such that $K_0(\cA)$ is isomorphic 
to the $\Gamma_2$-module $G'=(\Z/4\Z)^3$ in Example~\ref{Ex3}.
\end{example}

\begin{example}\label{Ex4}
Let $\Gamma_2\cong (\Z/2\Z)\times (\Z/2\Z)$ 
be as in the two previous examples. 
Let $G=\Z/8\Z$ be a $\Gamma_2$-module 
where an action of $\Gamma_2$ on $G$ is defined by 
$\sigma(g)=-g$ and $\tau(g)=3g$. 
In \cite[Example~2.7]{Ka5}, 
we saw that the $\Gamma_2$-module $G$ 
has no countable permutation presentations. 
Thus we cannot apply Theorem~\ref{Thm:constaction} 
nor Theorem~\ref{Thm:Main}. 
However Lemma~\ref{Lem:Z/2ZxGamma'} below 
together with Theorem~\ref{Thm:Main} 
shows that 
there exists a unital Kirchberg algebra $\cA$ 
in the Cuntz standard form 
and an action $\Gamma_2\act\cA$ such that 
$K_0(\cA)\cong G$ as $\Gamma_2$-modules and $K_1(\cA)=0$. 
By Theorem~\ref{Thm:KP}, 
we see that $\cA\cong \cOs_9$. 
This example shows that neither conditions in Theorem~\ref{Thm:constaction} 
nor Theorem~\ref{Thm:Main} is necessary. 
\end{example}

\begin{lemma}\label{Lem:Z/2ZxGamma'}
Let $\Gamma$ be a finite group, 
and $\cA$ be a unital Kirchberg algebra in the Cuntz standard form. 
Suppose that 
$\Gamma$ decomposes as a product $\Gamma=\Gamma_1\times \Gamma'$ 
where $\Gamma_1$ and $\Gamma'$ are subgroups of $\Gamma$ 
such that $\Gamma_1=\{1,\sigma\}\cong \Z/2\Z$. 
Then an action $\Gamma\act K_*(\cA)$ 
lifts to an action $\Gamma\act \cA$ 
if the restricted action $\Gamma'\act K_*(\cA)$ 
lifts to an action $\Gamma'\act \cA$ 
and $\sigma(x)=-x$ for all $x\in K_i(\cA)$ with $i=0,1$. 
\end{lemma}

\begin{proof}
By Theorem~\ref{Thm:Main}, 
there exists an action $\Gamma_1\act \cOs_\infty$ 
such that the induced action $\Gamma_1\act K_*(\cOs_\infty)$ satisfies 
$\sigma(x)=-x$ for $x\in K_0(\cOs_\infty)\cong\Z$. 
Then the tensor product action 
$\Gamma=\Gamma_1\times \Gamma'\act \cOs_\infty\otimes \cA\cong \cA$ 
is a lifting of the given action $\Gamma\act K_*(\cA)$ 
whose restriction to $\Gamma'$ lifts to $\Gamma'\act \cA$ 
and which satisfies $\sigma(x)=-x$ for all $x\in K_i(\cA)$ with $i=0,1$. 
\end{proof}

\begin{remark}\label{Rem:cannot}
A similar statement as this lemma holds for stable Kirchberg algebras. 
However it seems that the argument above does not work 
for unital Kirchberg algebras which is not 
in the Cuntz standard form. 
See the next example. 
\end{remark}

\begin{example}\label{Ex5}
The author does not know how to construct 
a unital Kirchberg algebra $\cA$ and an action $\Gamma_2\act\cA$ 
such that $K_0(\cA)$ is isomorphic to $\Gamma_2$-module $G$ 
in Example~\ref{Ex4} and $[1_{\cA}]\in K_0(\cA)$ corresponds to $4\in G$ 
which is fixed by the action of $\Gamma_2$. 
Thus it seems that in practice we need to worry about 
the position $[1_{\cA}]\in K_0(\cA)$ 
when considering the lifting problem. 
\end{example}

We finish this section 
by considering the splitting problem stated 
in the introduction of this paper. 
The following is easy to see. 

\begin{proposition}\label{Prop:split1}
For a Kirchberg algebra $\cA$ 
which has the isomorphic $K$-groups as $\cO_\infty$, 
the surjection 
$\Aut(\cA)\to \Aut(K_*(\cA))$ splits. 
\end{proposition}

\begin{proof}
Such a Kirchberg algebra $\cA$ is isomorphic to 
$M_k(\cO_\infty)$ for a positive integer $k$, 
$\cOs_\infty$ or $\cO_\infty\otimes\cK$ 
where $\cK$ is the \Ca of all compact operators on $\ell^2(\N)$. 
For $\cA=M_k(\cO_\infty)$, 
we have $K_*(\cA)=(\Z,k,0)$. 
Hence $\Aut(K_*(\cA))$ is trivial, 
and we need to do nothing. 
For $\cA=\cOs_\infty$ or $\cA=\cO_\infty\otimes\cK$, 
we have $\Aut(K_*(\cA))\cong \Z/2\Z$. 
Hence by Theorem~\ref{Thm:Main} 
the surjection 
$\Aut(\cA)\to \Aut(K_*(\cA))$ splits. 
\end{proof}

Take a positive integer $n$. 
For the Cuntz algebra $\cO_{n+1}$, 
the group $\Aut(K_*(\cO_{n+1}))$ is trivial. 
Thus the splitting problem for the Cuntz algebra $\cO_{n+1}$ 
are trivially valid. 
For the Kirchberg algebra $\cOs_{n+1}$ 
in the Cuntz standard form, 
we have $\Aut(K_*(\cOs_{n+1}))\cong \Aut(\Z/n\Z)$. 
The group $\Aut(\Z/n\Z)$ is isomorphic to 
the multiplicative group $(\Z/n\Z)^{\times}$ 
which can be computed as follows 
(see \cite[5.7.11 and 5.7.12]{Sc} for the proof). 
Let us write $n=2^{l_0}p_1^{l_1}p_2^{l_2}\cdots p_m^{l_m}$ 
for a non-negative integer $l_0$, positive integers $l_1,\ldots,l_m$ 
and distinct odd prime numbers $p_1,\ldots,p_m$. 
When $l_0=0$ or $1$ we have 
\[
(\Z/n\Z)^{\times} \cong 
(\Z/p_1^{l_1-1}(p_1-1)\Z)\times \cdots \times (\Z/p_m^{l_m-1}(p_m-1)\Z), 
\]
and when $l_0\geq 2$ we have 
\[
(\Z/n\Z)^{\times} \cong (\Z/2\Z)\times (\Z/2^{l_0-2}\Z)\times 
(\Z/p_1^{l_1-1}(p_1-1)\Z)\times \cdots \times (\Z/p_m^{l_m-1}(p_m-1)\Z). 
\]
Under these isomorphisms, 
the element $-1\in (\Z/n\Z)^{\times}$ corresponds to 
\[
\big(p_1^{l_1-1}(p_1-1)/2,\ldots, p_m^{l_m-1}(p_m-1)/2\big)
\] 
in the former case, 
and to 
\[
\big(1,0,p_1^{l_1-1}(p_1-1)/2,\ldots, p_m^{l_m-1}(p_m-1)/2\big)
\] 
in the latter case. 
From these computations, 
we get the following two lemmas whose proofs are omitted. 

\begin{lemma}\label{Lem:(Z/nZ)^x=cyc}
The group $(\Z/n\Z)^{\times}$ is cyclic 
if and only if 
\begin{itemize}
\item $n=1,2,4$ or 
\item 
$n=p^l,2p^l$ for an odd prime number $p$ and a positive integer $l$. 
\end{itemize}
\end{lemma}

\newpage

\begin{lemma}\label{Lem:(Z/nZ)^x=(Z/2Z)x(Z/2mZ)}
We have $(\Z/n\Z)^{\times}\cong (\Z/2\Z)\times (\Z/2m\Z)$ 
for a positive integer $m$ if and only if 
\begin{itemize}
\item $n=2^l$ for $l=3,4,\ldots$ or 
\item 
$n=4p^l$ for an odd prime number $p$ and a positive integer $l$ or  
\item $n=p_1^{l_1}p_2^{l_2}, 2p_1^{l_1}p_2^{l_2}$ 
for distinct odd prime numbers $p_1,p_2$ 
and positive integers $l_1,l_2$ satisfying 
\begin{itemize}
\item[$\cdot$] 
$(p_1-1)/2$ and $(p_2-1)/2$ are relatively prime and
\item[$\cdot$]
$p_1$ and $p_2-1$ are relatively prime if $l_1\geq 2$ and 
\item[$\cdot$]
$p_1-1$ and $p_2$ are relatively prime if $l_2\geq 2$.
\end{itemize}
\end{itemize}
In this case, 
we can take $-1\in (\Z/n\Z)^{\times}$ as 
the generator of the summand $\Z/2\Z$. 
\end{lemma}

On the splitting problem for the Kirchberg algebra $\cOs_{n+1}$, 
we get the following. 

\begin{proposition}\label{Prop:split2}
The surjection 
$\Aut(\cOs_{n+1})\to \Aut(K_*(\cOs_{n+1}))$ 
splits 
if $(\Z/n\Z)^{\times}$ is a cyclic group 
or a product of $\Z/2\Z$ and a cyclic group, i.e.  
\begin{itemize}
\item $n=2^l$ for $l\in\N$ or 
\item 
$n=p^l,2p^l,4p^l$ for an odd prime number $p$ and a positive integer $l$ or  
\item $n=p_1^{l_1}p_2^{l_2}, 2p_1^{l_1}p_2^{l_2}$ 
for $p_1,p_2,l_1,l_2$ as in Lemma~\ref{Lem:(Z/nZ)^x=(Z/2Z)x(Z/2mZ)}. 
\end{itemize}
\end{proposition}

\begin{proof}
Follows from Theorem~\ref{Thm:Main}, 
Lemma~\ref{Lem:Z/2ZxGamma'} 
Lemma~\ref{Lem:(Z/nZ)^x=cyc} and 
Lemma~\ref{Lem:(Z/nZ)^x=(Z/2Z)x(Z/2mZ)}. 
\end{proof}

Note that a finite abelian group satisfies 
the assumption of Theorem~\ref{Thm:Main} 
if and only if it is cyclic. 
Positive integers $n$ which do not 
satisfy the assumption in Proposition~\ref{Prop:split2} are 
\[
n= 24, 40, 48, 56, 60, 63, 65, 72, 80, 84, 85, 
   88, 91, 96, 104, 105, 112, 117, 120, \ldots. 
\]
It is clear that 
a similar statement of Proposition~\ref{Prop:split2} holds 
for stable Kirchberg algebras $\cO_{n+1}\otimes\cK$. 
For a unital Kirchberg algebras $M_k(\cO_{n+1})$ 
which is not in the Cuntz standard form, 
the author does not know the answer of 
the splitting problem 
except the case $\Aut\big(K_*(M_k(\cO_{n+1}))\big)$ is cyclic 
in which Theorem~\ref{Thm:Main} gives an affirmative answer. 
In particular, 
he does not know the answer for $M_4(\cO_{9})$ 
(Example~\ref{Ex5}). 
However the following proposition says that 
in generic cases $\Aut\big(K_*(M_k(\cO_{n+1}))\big)$ is cyclic 
and hence the splitting problem has an affirmative answer 
by Theorem~\ref{Thm:Main}. 
We omit the routine proof. 

\begin{proposition}\label{Prop:split3}
Let $n$ and $k$ be positive integers. 
Let us write $n=2^{l_0}p_1^{l_1}p_2^{l_2}\cdots p_m^{l_m}$ 
for $l_0\in\N$, positive integers $l_1,\ldots,l_m$ 
and distinct odd prime numbers $p_1,\ldots,p_m$. 
Suppose that $k$ is not divided by $p_i^{l_i}$ for $i=1,2,\ldots,m$, 
and not divided by $2^{l_0-1}$ when $l_0\geq 3$. 
Then 
\[
\Aut\big(K_*(M_k(\cO_{n+1}))\big)
=\Aut\big((\Z/n\Z,k,0)\big)
\cong \{x\in (\Z/n\Z)^{\times}\mid xk=k\}
\]
is a cyclic group, 
and hence the surjection 
$\Aut(M_k(\cO_{n+1}))\to \Aut\big(K_*(M_k(\cO_{n+1}))\big)$ 
splits by Theorem~\ref{Thm:Main}. 
\end{proposition}

\begin{remark}\label{Rem:Izumi}
In \cite[Subsection~6.3]{I}, 
Izumi discussed a problem 
which is similar to the splitting problem in this paper. 
His problem is equivalent to our splitting problem 
for a Kirchberg algebra $\cA$ with $K_i(\cA)=0$ 
for $i=0$ or $1$. 
By Proposition~\ref{Prop:split2}, 
the surjection 
$\Aut(\cOs_{n+1})\to \Aut(K_*(\cOs_{n+1}))$ 
splits if $n$ is a power of a prime number. 
This extends the latter statement of \cite[Theorem~6.12]{I} 
and the result explained in \cite[Remark~6.13]{I}. 
The former statement of \cite[Theorem~6.12]{I} can also be extended 
by the following corollary of Theorem~\ref{Thm:constaction}. 
\end{remark}

\begin{corollary}\label{Cor:Izumi}
Let $\Gamma$ be a finite group, and $\cA$ be a Kirchberg algebra. 
An action $\Gamma\act K_*(\cA)$ 
lifts to an action $\Gamma\act \cA$ 
if the induced two $\Gamma$-modules $K_0(\cA)$ and $K_1(\cA)$ are 
cohomologically trivial. 
\end{corollary}

\begin{proof}
By \cite[Theorem~VI.8.12]{Br} and \cite[Proposition~3.1]{Ka5}, 
every countable cohomologically trivial module has 
a countable permutation presentation. 
Hence Theorem~\ref{Thm:constaction} shows the conclusion. 
\end{proof}

\section{Relations to other classes of $C^*$-algebras}\label{Sec:Rel}

\subsection{Cuntz-Krieger algebras}\label{ssec:CK}

In this subsection, 
we see relations between our \CA s $\cO_{A,B}$ 
and the Cuntz-Krieger algebras $\cO_A$ 
introduced in \cite{CK} and generalized 
in \cite{EL,KPR} and other papers. 
Recall that the Cuntz-Krieger algebra $\cO_A$ 
of a matrix $A\in M_N(\N)$ is the universal \Ca 
generated by mutually orthogonal projections $\{p_i\}_{i=1}^N$ 
and partial isometries 
$\{\s{n}{i,j}\}_{(i,j)\in \Omega_{A},n\in\{1,2,\ldots,A_{i,j}\}}$ 
satisfying the relations (ii) and (iii) in Definition~\ref{Def:OAB}. 
When $N=\infty$, $\cO_A$ should be called the Exel-Laca algebra 
or the graph algebra of the row-finite graph associated with $A$ 
(see \cite{EL} or \cite[Chapter~2]{Ra}). 
Take $N\in\fC$ and 
$A,B\in M_N(\Z)$ satisfying the condition (0), 
and fix them. 
We will show that the natural map $\cO_A\to \cO_{A,B}$ is injective. 
To this end, 
we need the following definition and proposition. 

\begin{definition}\label{Def:gauge}
For each $z\in\T$ 
the universality of $\cO_{A,B}$ 
shows the existence of an automorphism $\beta_z$ of $\cO_{A,B}$ 
such that $\beta_z(p_i)=p_i$, $\beta_z(u_i)=u_i$ for each $i$ 
and $\beta_z(\s{n}{i,j})=z\s{n}{i,j}$ 
for each $(i,j)\in \Omega_{A}$ and $n\in\Z$. 
We call the action $\beta\colon \T\act \cO_{A,B}$ 
the {\em gauge action} of $\cO_{A,B}$. 
\end{definition}

\begin{proposition}\label{Prop:repOAB}
For each $i\in \{1,2,\ldots,N\}$, 
the projection $p_i\in \cO_{A,B}$ is non-zero. 
\end{proposition}

\begin{proof}
It suffices to find operators 
$\{P_i\}_{i=1}^N$, $\{U_i\}_{i=1}^N$, 
and $\{\bS{n}{i,j}\}_{(i,j)\in \Omega_{A},n\in\Z}$ 
on some Hilbert space 
such that they satisfy the relations (i) to (iii) in Definition~\ref{Def:OAB}
and $P_i\neq 0$ for all $i$. 

We set 
\[
\Omega_{A}^{\infty}:=\big\{\mu=(\mu_0,\mu_1,\ldots)\in \{1,2,\ldots,N\}^\infty
\ \big|\ \text{$(\mu_{l-1},\mu_{l})\in \Omega_{A}$ for $l=1,2,\ldots$}\big\}. 
\]
For each $\mu\in \Omega_{A}^{\infty}$, 
we define a group $F_\mu$ to be the abelian group 
generated by generators $\{\xi_{\mu}^{l}\}_{l=1}^\infty$ 
with relations 
$A_{\mu_{l-1},\mu_{l}}\xi_{\mu}^{l}=B_{\mu_{l},\mu_{l+1}}\xi_{\mu}^{l+1}$ 
for $l=1,2,\ldots$. 
For $(i,j)\in \Omega_{A}$ 
and $\mu\in \Omega_{A}^{\infty}$ with $\mu_0=j$, 
we set $i\mu\in \Omega_{A}^{\infty}$ by $(i\mu)_0=i$ 
and $(i\mu)_l=\mu_{l-1}$ for $l=1,2,\ldots$. 
We define a homomorphism $\sigma_{i}\colon F_\mu\to F_{i\mu}$
by $\sigma_{i}(\xi_{\mu}^{l})=\xi_{i\mu}^{l+1}$ for all $l=1,2,\ldots$. 
It is easy to see that $\sigma_{i}$ is a well-defined injective 
homomorphism. 
It is also easy to see that 
we have $\sigma_{i}(B_{\mu_0,\mu_1}\xi_{\mu}^{1})=A_{i,j}\xi_{i\mu}^{1}$, 
and the quotient $F_{i\mu}/\sigma_{i}(F_\mu)$ 
is isomorphic to $\Z/A_{i,j}\Z$ which is generated 
by the image of $\xi_{i\mu}^{1}\in F_{i\mu}$. 
Hence we have 
\[
F_{i\mu}=\coprod_{n=1}^{A_{i,j}}\big(\sigma_{i}(F_\mu)+n\xi_{i\mu}^{1}\big).
\]
Let $\F^\infty:=\coprod_{\mu\in \Omega_{A}^{\infty}}F_\mu$. 
We define operators $\{P_i\}_{i=1}^N$, $\{U_i\}_{i=1}^N$, 
and $\{\bS{n}{i,j}\}_{(i,j)\in \Omega_{A},n\in\Z}$ 
on the Hilbert space $\ell^2(\F^\infty)$ by 
\[
P_i(\delta_{\xi}):=
\begin{cases}
\delta_{\xi} & \text{if $\mu_0=i$}\\
0 & \text{if $\mu_0\neq i$},
\end{cases}\qquad 
U_i(\delta_{\xi}):=
\begin{cases}
\delta_{\xi+B_{\mu_0,\mu_1}\xi_{\mu}^{1}} & \text{if $\mu_0=i$}\\
0 & \text{if $\mu_0\neq i$},
\end{cases}
\]
\[
\bS{n}{i,j}(\delta_{\xi}):=
\begin{cases}
\delta_{\sigma_{i}(\xi)+n\xi_{i\mu}^{1}} & \text{if $\mu_0=j$}\\
0 & \text{if $\mu_0\neq j$},
\end{cases}
\]
for $\xi\in F_\mu\subset\F^\infty$. 
Then it is routine to check that these operators satisfy 
the relations (i) to (iii) in Definition~\ref{Def:OAB}. 
For each $i\in \{1,2,\ldots,N\}$ 
we can find $\mu\in \Omega_{A}^{\infty}$ with $\mu_0=i$ 
because $\Omega_A(i)\neq\emptyset$ for all $i$. 
Hence $P_i\neq 0$ for all $i$. 
We are done. 
\end{proof}

As a corollary of this proposition, 
we get the following. 
We say that 
$B\in M_N(\Z)$ {\em has a zero row} 
if there exists $i\in\{1,2,\ldots,N\}$ 
with $B_{i,j}=0$ for all $j\in\{1,2,\ldots,N\}$. 
Otherwise we say that $B\in M_N(\Z)$ {\em has no zero rows}. 

\begin{corollary}\label{Cor:spec(u_i)}
The spectrum of $u_i\in \cO_{A,B}$ contains $\T$ 
for all $i\in\{1,2,\ldots,N\}$ 
if and only if $B$ has no zero rows. 
\end{corollary}

\begin{proof}
If there exists $i\in\{1,2,\ldots,N\}$ 
with $B_{i,j}=0$ for all $j\in\{1,2,\ldots,N\}$, 
then we have 
\[
u_i=u_ip_i=u_i\Big(\sum_{j\in \Omega_{A}(i)}
\sum_{n=1}^{A_{i,j}}\s{n}{i,j}\s{n}{i,j}^*\Big)
=\sum_{j\in \Omega_{A}(i)}\sum_{n=1}^{A_{i,j}}\s{n}{i,j}\s{n}{i,j}^*
=p_i. 
\]
Hence if the spectrum of $u_i\in \cO_{A,B}$ contains $\T$ 
for all $i$, 
then $B$ has no zero rows. 
Conversely suppose that $B$ has no zero rows. 
Take $i\in\{1,2,\ldots,N\}$. 
To show that the spectrum of $u_i\in \cO_{A,B}$ contains $\T$, 
it suffices to show that 
the spectrum of the operator $U_i$ on $\ell^2(\F^\infty)$ 
defined in the proof of Proposition~\ref{Prop:repOAB} 
contains $\T$. 
Since $B$ has no zero rows, 
we can find $\mu\in \Omega_{A}^{\infty}$ 
with $\mu_0=i$ and $B_{\mu_{l-1},\mu_{l}}\neq 0$ for all $l$. 
Then the order of the element $\xi:=B_{\mu_0,\mu_1}\xi_{\mu}^{1}$ 
in the group $F_\mu$ is infinite. 
Hence the spectrum of $U_i$ contains $\T$ 
because $U_i$ acts as a bilateral shift on 
$\ell^2(\Z\xi)\cong\ell^2(\Z)$. 
This completes the proof. 
\end{proof}

\begin{remark}
From the operators on $H:=\ell^2(\F^\infty)$ 
defined in the proof of Proposition~\ref{Prop:repOAB}, 
we get a representation $\pi\colon \cO_{A,B}\to B(H)$. 
In general, 
the representation $\pi$ is not faithful. 
However 
by computing the spectrum of $u_i$ and $U_i$, 
we can show that 
the restriction of $\pi$ to the \Ca 
generated by $\{u_i\}_{i=1}^N$ is injective. 
From this fact and a variant of gauge-invariant uniqueness theorems, 
we can show that the \shom 
\[
\widetilde{\pi}\colon \cO_{A,B}\to C\big(\T,B(H)\big)
\cong C(\T)\otimes B(H) \subset B\big(L^2(\T)\otimes H\big)
\] 
defined by $\widetilde{\pi}(x)(z)=\pi\big(\beta_z(x)\big)$ 
is injective. 
Thus we get an explicit faithful representation of $\cO_{A,B}$. 
Note that 
a variant of Cuntz-Krieger uniqueness theorems 
shows that the representation $\pi$ itself is faithful 
when $A$ and $B$ are sufficiently complicated. 
\end{remark}

\begin{proposition}\label{Prop:OA->OA}
The natural \shom $\cO_A\to \cO_{A,B}$ defined by 
$p_i\mapsto p_i$ and $\s{n}{i,j}\mapsto \s{n}{i,j}$ 
is injective. 
\end{proposition}

\begin{proof}
This follows from 
the gauge-invariant uniqueness theorem of the Cuntz-Krieger algebras 
(see \cite[Theorem~2.2]{Ra}) 
with the help of the gauge action $\beta$ of $\cO_{A,B}$ 
and Proposition~\ref{Prop:repOAB}. 
\end{proof}

By this proposition, 
we can consider the Cuntz-Krieger algebra $\cO_A$ 
as a \Csa of $\cO_{A,B}$. 

It is routine to see that 
$\sum_{i=1}^N p_i$ converges to the unit of 
the multiplier algebra $\cM(\cO_{A,B})$ of $\cO_{A,B}$ 
in the strict topology 
(when $N<\infty$, 
$\sum_{i=1}^N p_i$ is the unit of $\cO_{A,B}$). 
Hence $\sum_{i=1}^N u_i$ 
converges to a unitary $u$ of $\cM(\cO_{A,B})$ 
in the strict topology. 
The unitary $u\in \cM(\cO_{A,B})$ 
has the following commutation relations 
with the generators of $\cO_A\subset \cO_{A,B}$; 
\[
up_i=p_iu\ (=u_i),\qquad 
u^k\s{n}{i,j}=\s{n'}{i,j}u^{k'}\ (=\s{m}{i,j})
\]
for all $i\in\{1,2,\ldots,N\}$, 
$(i,j)\in\Omega_{A}$, $n,n'\in\{1,2,\ldots,A_{i,j}\}$ 
and $k,k'\in\Z$ 
with $n+kB_{i,j}=n'+k'A_{i,j}\ (=m)$. 
The following proposition says that 
the \Ca $\cO_{A,B}$ can be defined to be 
the universal \Ca generated 
by the product of the Cuntz-Krieger $\cO_A\subset \cO_{A,B}$
and the unitary $u\in \cM(\cO_{A,B})$ 
with the commutation relations above. 

\begin{proposition}\label{Prop:OA_u}
For a non-degenerate representation $\cO_A\to B(H)$ 
and a unitary $U\in B(H)$ satisfying 
the commutation relations above, 
there exists a unique representation $\cO_{A,B}\to B(H)$ 
extending the representation $\cO_A\to B(H)$ 
and sending the unitary $u\in \cM(\cO_{A,B})$ 
to $U\in B(H)$. 
\end{proposition}

\begin{proof}
Straightforward. 
\end{proof}

When $B=0$ or $B=A$, 
the \Ca $\cO_{A,B}$ can be explicitly described 
using the Cuntz-Krieger algebra $\cO_A$. 

\begin{proposition}\label{Prop:A0AA}
For a matrix $A\in M_N(\N)$, 
we have $\cO_{A,0}\cong \cO_A$ and 
$\cO_{A,A}\cong \cO_A\otimes C(\T)$. 
\end{proposition}

\begin{proof}
When $B=0$, 
we have $u_i=p_i$ for all $i$ 
as computed in the proof of Corollary~\ref{Cor:spec(u_i)}. 
This shows that the \shom $\cO_A\to \cO_{A,0}$ 
in Proposition~\ref{Prop:OA->OA} is surjective, 
and hence an isomorphism. 
Note that in this case we have $u=1$. 

When $B=A$, 
the unitary $u\in \cM(\cO_{A,A})$ commutes with $\cO_A\subset \cO_{A,A}$. 
Hence $\cO_{A,A}\cong \cO_A\otimes C(\T)$ by Proposition~\ref{Prop:OA_u}. 
\end{proof}

\subsection{Topological graph algebras}\label{ssec:TopG}

In \cite{Ka1}, 
a notion of topological graphs was introduced, 
and a construction of a \Ca $\cO(E)$ from a topological graph $E$ 
was given. 
This construction generalizes 
the one of Cuntz-Krieger algebras or 
more generally of graph algebras. 
In \cite[Definition~6.1]{Ka4} 
a topological graph $E_{A,B}$ was constructed 
from two matrices $A\in M_\infty(\N)$ and $B\in M_\infty(\Z)$ 
such that $A_{i,j}=0$ implies $B_{i,j}=0$, 
and in \cite[Proposition~B.1]{Ka4} 
generators and relations of the \Ca $\cO(E_{A,B})$ 
were provided. 
The definition of the \Ca $\cO_{A,B}$ 
is motivated by these generators and relations 
of $\cO(E_{A,B})$, 
and in fact by \cite[Proposition~B.2]{Ka4} 
we have $\cO_{A,B}\cong \cO(E_{A,B})$ 
for $A,B\in M_\infty(\Z)$ satisfying the condition (0) 
if $B$ has no zero rows. 
As pointed out in \cite[Remark~6.3]{Ka4}, 
a similar construction of a topological graph $E_{A,B}$ 
from $A,B\in M_N(\Z)$ satisfying the condition (0)
for $N<\infty$ is possible, 
and we can show $\cO_{A,B}\cong \cO(E_{A,B})$ 
if $B$ has no zero rows. 
However in the case that $B$ has a zero row, 
the natural surjection from $\cO(E_{A,B})$ to $\cO_{A,B}$ 
is never injective 
because there exists a restriction ``$v\in E^0_{\text{rg}}\cap E^0_m$'' 
in the relation (iii) of \cite[Proposition~B.2]{Ka4} 
although we do not consider the corresponding restriction 
when defining $\cO_{A,B}$. 
This restriction looks natural from the graph algebraic point of view 
(cf.\ the fact that the Cuntz-Krieger relation (CK2) is not assumed 
at sources in p.6 of \cite{Ra}). 
We do not adopt this restriction 
because this seems to be unnatural 
from our point of view in this paper, 
and makes the computation of $K$-groups 
in Proposition~\ref{Prop:KOAB} complicated. 

In the case that $B$ has no zero rows, 
the computation of $K$-groups 
in Proposition~\ref{Prop:KOAB} 
without considering $\Gamma_{A,B}$-actions 
follows from \cite[Lemma~6.2]{Ka4}, 
and Proposition~\ref{Prop:OAB_sn} 
follows from \cite[Propositions~6.1 and 6.6]{Ka1}. 
Even when $B$ has a zero row, 
we can get Proposition~\ref{Prop:OAB_sn} 
using the theory of topological graph algebras 
with extra efforts. 
We give direct complete proofs of them in the following three sections 
without this restriction on $B$ 
because we naturally get it on the way to other results. 

We can show Proposition~\ref{Prop:OAB_Ki} 
from \cite[Proposition~6.5]{Ka4} if $N=\infty$, 
and by repeating a similar argument 
we can show it for $N<\infty$. 
For the readers' convenience, 
we give a self-contained proof of Proposition~\ref{Prop:OAB_Ki} 
in Section~\ref{Sec:pi} 
instead of forcing the readers follow long arguments 
from \cite{Ka1} to \cite{Ka4} 
which are too long and too general for our purpose.

\subsection{Cuntz-Pimsner algebras}\label{ssec:CP}

In \cite{Pi}, 
Pimsner introduced a construction of a \Ca $\cO_\sX$ 
from a \Cc $\sX$ which was called a Hilbert bimodule in \cite{Pi}. 
This \Ca $\cO_\sX$ is now called a {\em Cuntz-Pimsner algebra}, 
and it generalizes Cuntz-Krieger algebras. 
From $N\in\fC$ and 
two matrices $A,B\in M_N(\Z)$ satisfying the condition (0), 
we can define a full \Cc 
$\sX_{A,B}\,\big({\cong\bigoplus_{(i,j)\in\Omega_A} C(\T)}\big)$ 
over a commutative \Ca $\cA_N\cong C_0(\{1,2,\ldots,N\}\times \T)$
in a similar way as in \cite{Ka4}, 
and show that the Cuntz-Pimsner algebra $\cO_{\sX_{A,B}}$ 
is isomorphic to our \Ca $\cO_{A,B}$. 
We remark that in \cite{Ka2} 
a modified construction of Cuntz-Pimsner algebras was proposed 
so that it generalizes topological graph algebras, 
but here we mean the original definition of Cuntz-Pimsner algebras. 
Hence the natural \shom $\cA_N\to \cO_{\sX_{A,B}}$ is injective 
if and only if the left action of the \Cc $\sX_{A,B}$ is faithful 
which turns out to be equivalent that $B$ has no zero rows 
(cf.\ Corollary~\ref{Cor:spec(u_i)}). 
We also note that the Toeplitz algebra $\cT_{\sX_{A,B}}$ 
defined in \cite{Pi} is isomorphic to the \Ca $\cT_{A,B}$ 
which will be defined in Definition~\ref{Def:TAB}, 
and that we have 
\begin{align*}
\cA_N&\cong
\cspa\big\{u_i^n\mid i\in\{1,2,\ldots,N\},n\in\Z\big\} \subset \cT_{A,B}, \\
\sX_{A,B}&\cong
\cspa\{\s{n}{i,j}\mid (i,j)\in \Omega_{A},n\in\Z\} \subset \cT_{A,B}. 
\end{align*}
In the proof of Proposition~\ref{Prop:KOAB}, 
we use the brilliant idea of Pimsner in \cite[Section~4]{Pi} 
for computing the $K$-groups of Cuntz-Pimsner algebras. 
We repeat Pimsner's argument in our terminology in Section~\ref{Sec:K} 
because the left actions of our \Cc $\sX_{A,B}$ 
need not be faithful 
although it was assumed to be faithful in \cite{Pi}, 
and because 
we also have to compute $\Gamma_{A,B}$-actions on 
the $K$-groups of $\cO_{A,B}$. 
When $B$ has no zero rows, 
Proposition~\ref{Prop:OAB_sn} follows from 
the theory of Cuntz-Pimsner algebras (see \cite{Ka3} for example).

\section{Structures and the nuclearity of $\cO_{A,B}$}\label{Sec:nuc}

Let $A,B\in M_N(\Z)$ satisfy the condition (0). 
In this section, 
we define a dense $*$-algebra $\acO_{A,B}$ of $\cO_{A,B}$ 
whose elements are described explicitly, 
and show that $\cO_{A,B}$ is nuclear 
by examining the so-called core $\cO_{A,B}^{\T}$ of $\cO_{A,B}$. 
The analysis in this section will be used in the next section. 

We set $\Omega_{A}^0:=\{\mu=(\mu_0)\mid \mu_0\in \{1,2,\ldots,N\}\}$, 
and 
\[
\Omega_{A}^{k}:=\{\mu=(\mu_0,\mu_1,\ldots,\mu_k)\in \{1,2,\ldots,N\}^{k+1}
\mid \text{$(\mu_{l-1},\mu_{l})\in \Omega_{A}$ for $l=1,2,\ldots,k$}\}, 
\]
for a positive integer $k$. 
We define $\Omega_{A}^{*}:=\coprod_{k\in\N}\Omega_{A}^{k}$. 

For $\mu\in \Omega_{A}^0$, 
let $F_{\mu}\cong \Z$ be the infinite cyclic group 
whose generator is given by $\xi_{\mu}^{0}\in F_{\mu}$. 
For $\xi=n\xi_{\mu}^{0}\in F_{\mu}$ with $n\in\Z$, 
we define $s_{\xi}:=u_{\mu_0}^n\in\cO_{A,B}$. 
Let us take $\mu\in \Omega_{A}^{k}$ 
for a positive integer $k$. 
Let $F_\mu$ be the abelian group generated 
by elements $\{\xi_{\mu}^{l}\}_{l=1}^{k}$ with relations 
$A_{\mu_{l-1},\mu_{l}}\xi_{\mu}^{l}=B_{\mu_{l},\mu_{l+1}}\xi_{\mu}^{l+1}$ 
for $l=1,2,\ldots,k-1$. 
For $\xi=\sum_{l=1}^kn_l\xi_{\mu}^{l}\in F_\mu$ 
with $n_1,n_2,\ldots,n_k\in\Z$, 
we define $s_{\xi}\in\cO_{A,B}$ by 
\[
s_{\xi}:=\s{n_1}{\mu_0,\mu_1}\s{n_2}{\mu_1,\mu_2}\cdots 
\s{n_k}{\mu_{k-1},\mu_{k}}. 
\]
It is routine to check that this is well-defined, 
and $s_{\xi}$ satisfies $s_{\xi}^*s_{\xi}=p_{\mu_k}$. 

We set 
$\F^k:=\coprod_{\mu\in \Omega_{A}^{k}}F_\mu$ for $k\in\N$ 
and $\F^*:=\coprod_{\mu\in \Omega_{A}^{*}}F_\mu$. 
Thus $\F^*=\coprod_{k\in\N}\F^k$. 
We define 
\[
\acO_{A,B}
:=\spa\big\{s_{\xi}s_{\eta}^*\ \big|\ \xi,\eta\in \F^{*}\}. 
\]
By the following lemma, 
$\acO_{A,B}$ is the $*$-algebra generated 
by the generators $\{p_i,u_i,\s{n}{i,j}\}$, 
and hence dense in $\cO_{A,B}$. 

\begin{lemma}\label{Lem:sxsx*}
For $\xi\in \F^k$ and $\eta\in\F^l$ with $k,l\in\N$, 
we have the following. 
\begin{itemize}
\item We have either $s_{\xi}s_{\eta}=0$ or 
$s_{\xi}s_{\eta}=s_{\zeta}$ for some $\zeta\in \F^{k+l}$. 
\item When $k\geq l$, 
we have either $s_{\eta}^*s_{\xi}=0$ or 
$s_{\eta}^*s_{\xi}=s_{\zeta}$ for some $\zeta\in \F^{k-l}$. 
\item 
When $k=l\geq 1$, 
we have $s_{\eta}^*s_{\xi}\neq 0$ if and only if 
$\xi,\eta\in F_\mu$ 
and $\xi-\eta=nA_{\mu_{k-1},\mu_k}\xi_{\mu}^{k}$ for $n\in\Z$. 
In this case, we have $s_{\eta}^*s_{\xi}=u_{\mu_k}^{n}$. 
\end{itemize} 
\end{lemma}

\begin{proof}
Straightforward. 
\end{proof}

Let $\beta\colon \T\act \cO_{A,B}$ be 
the gauge action of $\cO_{A,B}$ 
defined in Definition~\ref{Def:gauge}. 
We examine the structure of 
the fixed point algebra $\cO_{A,B}^{\T}$ of the gauge action $\beta$. 
For $k\in\N$, 
we define $\Cr_{A,B}^k\subset \cO_{A,B}$ by 
\[
\Cr_{A,B}^k:=\cspa\{s_{\xi}s_{\eta}^*\mid \xi,\eta\in \F^k\}. 
\]
By Lemma~\ref{Lem:sxsx*}, 
$\Cr_{A,B}^k$ is a \Csa of $\cO_{A,B}$.

\begin{lemma}\label{Lem:Ck}
For $k\in\N$, 
the \Ca $\Cr_{A,B}^k$ is isomorphic to 
$\bigoplus_{i=1}^N C^*(u_i)\otimes \cK\big(\ell^2(X^k_i)\big)$ 
where $\{X^k_i\}_{i=1}^N$ are subsets of $\F^k$ defined in the proof. 
\end{lemma}

\begin{proof}
For $k\in\N$ and $i\in\{1,2,\ldots,N\}$, 
we define $(\Omega_{A}^{k})_i\subset \Omega_{A}^{k}$ by 
\[
(\Omega_{A}^{k})_i
:=\{\mu\in \Omega_{A}^{k}\mid \mu_k=i\}.
\]
We define $\F^k_i\subset \F^k$ by 
$\F^k_i=\coprod_{\mu\in (\Omega_{A}^{k})_i}F_\mu$.
Note that we have 
$\F^k=\coprod_{i=1}^N \F^k_i$ and 
$\F^k_i=\{\xi\in \F^k\mid s_{\xi}^*s_{\xi}=p_i\}$ for each $i$. 
We define $(\Cr_{A,B}^k)_i\subset \Cr_{A,B}^k$ by 
\[
(\Cr_{A,B}^k)_i:=\cspa\{s_{\xi}s_{\eta}^*\mid \xi,\eta\in \F^k_i\}. 
\]
Then 
we have $\Cr_{A,B}^k\cong \bigoplus_{i=1}^N(\Cr_{A,B}^k)_i$. 

For $\mu\in \Omega_{A}^{0}$, 
we define $X_\mu:=\{0\}\subset F_\mu$. 
For $\mu\in \Omega_{A}^{k}$ 
with a positive integer $k$, 
let $X_\mu\subset F_\mu$ be the image of the map 
\[
\prod_{l=1}^k\{1,2,\ldots,A_{\mu_{l-1},\mu_{l}}\} \ni (n_l)_{l=1}^k 
\mapsto \sum_{l=1}^k n_l\xi_{\mu}^{l}\in F_\mu. 
\]
Then it is routine to see that the map above is injective, 
and each $\xi\in F_\mu$ can be uniquely written as 
$\xi=\eta+nA_{\mu_{k-1},\mu_k}\xi_{\mu}^{k}$ 
for $\eta\in X_\mu$ and $n\in\Z$. 
Note that we have $s_\xi=s_\eta u_{\mu_k}^n$. 
For $k\in\N$ and $i\in\{1,2,\ldots,N\}$, 
we define a subset $X^k_i\subset \F^k_i$ 
by $X^k_i:=\coprod_{\mu\in (\Omega_{A}^{k})_i}X_\mu$. 
Then for $\xi,\eta\in X^k_i$ with $\xi\neq\eta$ 
we have $s_{\eta}^*s_{\xi}=0$ by Lemma~\ref{Lem:sxsx*}. 
Now it is routine to check that the map 
\[
C^*(u_i)\otimes \cK\big(\ell^2(X^k_i)\big)
\ni x\otimes \theta_{\xi,\eta}
\mapsto s_{\xi}xs_{\eta}^*\in (\Cr_{A,B}^k)_i
\]
is an isomorphism 
where $\{\theta_{\xi,\eta}\}_{\xi,\eta\in X^k_i}$ 
are the matrix units of $\cK\big(\ell^2(X^k_i)\big)$. 
We are done. 
\end{proof}

We use the following lemma in the next section. 

\begin{lemma}\label{Lem:bab}
For a positive element $a\in \Cr_{A,B}^k$ with $k\in\N$, 
there exists $b\in \cspa\{s_{\xi}\mid \xi\in \F^k\}$ 
such that $\|b\|=1$, $\|b^*ab\|=\|a\|$ 
and $b^*ab\in C^*(u_i)$ 
for some $i\in\{1,2,\ldots,N\}$. 
\end{lemma}

\begin{proof}
We use the notation established 
in the proof of Lemma~\ref{Lem:Ck}. 

Take a positive element $a\in \Cr_{A,B}^k$ with $k\in\N$. 
There exist positive elements 
$a_i\in (\Cr_{A,B}^k)_i$ for $i\in\{1,2,\ldots,N\}$ 
such that $a=\sum_{i=1}^N a_i$ and $\|a\|=\max_{i}\|a_i\|$. 
Take $i$ with $\|a\|=\|a_i\|$. 
For each $z\in\T$ in the spectrum of $u_i$, 
we define a $*$-homomorphism 
$\varphi_z\colon(\Cr_{A,B}^k)_i\to \cK\big(\ell^2(X^k_i)\big)$ 
by $\varphi_z(s_{\xi}u_i^ns_{\eta}^*)=z^n\theta_{\xi,\eta}$ 
for $\xi,\eta\in X^k_i$ and $n\in\Z$. 
It is well-defined and 
we have $\|a\|=\|a_i\|=\max_{z}\|\varphi_z(a_i)\|$ 
by the proof of Lemma~\ref{Lem:Ck}. 
Take $z\in\T$ with $\|a\|=\|\varphi_z(a_i)\|$. 
Since $\varphi_z(a_i)\in \cK\big(\ell^2(X^k_i)\big)$ 
is positive, 
there exists $w\in \ell^2(X^k_i)$ 
such that $\|w\|=1$ and 
$\|a\|=\|\varphi_z(a_i)\|=\ip{w}{\varphi_z(a_i)w}$ 
where $\ip{\cdot}{\cdot}$ denotes 
the inner product in $\ell^2(X^k_i)$. 
Express $w\in \ell^2(X^k_i)$ as 
$w=\sum_{\xi\in X^k_i}\lambda_\xi\delta_\xi$ 
using $\lambda_\xi\in\C$ with $\sum_{\xi\in X^k_i}|\lambda_\xi|^2=1$. 
We set $b=\sum_{\xi\in X^k_i}\lambda_\xi s_\xi\in \cO_{A,B}$ 
which converges in the norm topology 
and satisfies $\|b\|=1$. 
By Lemma~\ref{Lem:sxsx*}, 
$b^*ab=b^*a_{i}b\in C^*(u_i)$. 
We obtain $\psi_z(b^*a_{i}b)=\ip{w}{\varphi_z(a_i)w}=\|a\|$ 
where $\psi_z\colon C^*(u_i)\to \C$ is defined by $\psi_z(u_i)=z$. 
Hence we get $\|b^*ab\|=\|a\|$. 
Thus the element $b\in \cO_{A,B}$ satisfies the desired conditions. 
\end{proof}

\begin{definition}\label{Def:condexp}
We define 
a faithful conditional expectation 
$\varPhi\colon \cO_{A,B}\to \cO_{A,B}^{\T}$ by 
\[
\varPhi(x):=\int_{\T}\beta_z(x)dz
\]
where $dz$ is the normalized Haar measure on $\T$. 
\end{definition}

For $\xi\in \F^k$, we get $\beta_z(s_{\xi})=z^ks_{\xi}$. 
Hence for $\xi\in \F^k$ and $\eta\in \F^l$ we have 
\[
\varPhi\big(s_{\xi}s_{\eta}^*\big)=
\begin{cases}
s_{\xi}s_{\eta}^* & \text{if $k=l$},\\
0 & \text{if $k\neq l$}.
\end{cases}
\]

\begin{lemma}\label{Lem:indlim}
We have 
$\Cr_{A,B}^{k}\subset \Cr_{A,B}^{k+1}$ for each $k\in\N$, 
and $\cO_{A,B}^{\T}=\overline{\bigcup_{k\in\N}\Cr_{A,B}^{k}}$.
\end{lemma}

\begin{proof}
By the relation (iii) in Definition~\ref{Def:OAB}, 
we have $\Cr_{A,B}^{k}\subset \Cr_{A,B}^{k+1}$. 
Since the conditional expectation $\varPhi$ is bounded, 
$\varPhi(\acO_{A,B})$ is dense in $\cO_{A,B}^{\T}$. 
By the computation of $\varPhi(s_{\xi}s_{\eta}^*)$ above, 
we get $\varPhi(\acO_{A,B})\subset \bigcup_{k\in\N}\Cr_{A,B}^{k}$. 
Thus we obtain $\cO_{A,B}^{\T}=\overline{\bigcup_{k\in\N}\Cr_{A,B}^{k}}$. 
\end{proof}

\begin{proposition}\label{Prop:OAB_nuc}
For $A,B\in M_N(\Z)$ satisfying the condition {\rm (0)}, 
the \Ca $\cO_{A,B}$ is nuclear. 
\end{proposition}

\begin{proof}
By Lemma~\ref{Lem:Ck} and Lemma~\ref{Lem:indlim}, 
the \Ca $\cO_{A,B}^{\T}$ is nuclear. 
Hence $\cO_{A,B}$ is nuclear by \cite[Proposition~2]{DLRZ}. 
\end{proof}

\section{Pure infiniteness of $\cO_{A,B}$}\label{Sec:pi}

In this section, 
we prove Proposition~\ref{Prop:OAB_Ki}. 
We need the following lemma. 

\begin{lemma}\label{Lem:isometry}
Let $K$ be an integer greater than $1$. 
Let $\cA$ be a unital \Ca 
which has a unitary $u\in \cA$ and an isometry $s\in \cA$ 
satisfying that $s^*u^ns=0$ for $n=1,2,\ldots,K-1$ and $u^Ks=su$. 
Then for $k\in\N$ and a function $f\in C(\T)$ which is $1$ 
on some non-empty open subset of $\T$, 
there exist $k_0\in\N$ with $k_0\geq k$ 
and a function $g\in C(\T)$ such that 
$v:=g(u)s^{k_0}\in \cA$ is an isometry satisfying $f(u)v=v$ 
and $v^*u^ns^lv=0$ for all $n\in\Z$ and $l\in\{1,2,\ldots,k\}$. 
\end{lemma}

\begin{proof}
We define a linear map $E\colon C(\T)\to C(\T)$ by 
\[
E(g)(z):=\frac{1}{K}\sum_{\text{$w\in\T$ with $w^K=z$}}g(w)
\]
for $g\in C(\T)$ and $z\in\T$. 
Then the map $E$ is characterized 
by the contracting linear map satisfying 
\[
E(z^n)=\begin{cases}
z^{m}& \text{if $n=Km$ for $m\in\Z$}\\
0 & \text{if $n\in \Z\setminus K\Z$,}
\end{cases}
\]
where $z\in C(\T)$ is the identity function. 
Since we have $s^*u^{Km}s=u^{m}$ for $m\in\Z$ 
and $s^*u^ns=0$ for $n\in \Z\setminus K\Z$, 
we get $s^*g(u)s=E(g)(u)$ for all $g\in C(\T)$. 

Take $k\in\N$ and a function $f\in C(\T)$ which is $1$ 
on some non-empty open subset of $\T$. 
Then we can find a real number $t_0$, 
and $k_0\in\N$ with $k_0\geq k$ 
such that $f$ is $1$ on 
\[
I_0:=\big\{e^{2\pi \sqrt{-1}t}\ \big|\ 
t_0-K^{-k_0}\leq t\leq t_0+K^{-k_0}\big\} 
\subset \T, 
\]
and $I_0\cap I_l=\emptyset$ for $l=1,2,\ldots,k$ where 
\[
I_l:=\big\{e^{2\pi \sqrt{-1}K^{l}t}\ \big|\ 
t_0-K^{-k_0}\leq t\leq t_0+K^{-k_0}\big\} 
\subset \T. 
\]
Let us define a positive function $g_0\in C(\T)$ by 
\[
g_0(z):=\begin{cases}
K^{k_0}-K^{2k_0}|t-t_0|& 
(\text{for $z=e^{2\pi \sqrt{-1}t}$ 
with $t_0-K^{-k_0}\leq t\leq t_0+K^{-k_0}$}),\\
0&(\text{for $z\notin I_0$}).
\end{cases}
\]
Then for all $z\in\T$ we have 
\[
E^{k_0}(g_0)(z)
=\frac{1}{K^{k_0}}\sum_{\text{$w\in\T$ with $w^{K^{k_0}}=z$}}g_0(w)
=1. 
\]
Thus we obtain $E^{k_0}(g_0)=1$. 
We set $g:=g_0^{1/2}\in C(\T)$.
Then the element $v:=g(u)s^{k_0}\in \cA$ 
satisfies $f(u)v=v$ because $fg=g$. 
From the computation $v^*v=(s^*)^{k_0}g_0(u)s^{k_0}=E^{k_0}(g_0)(u)=1$, 
we see that $v$ is an isometry. 
We will show $g(u)s^lg(u)=0$ for $l\in\{1,2,\ldots,k\}$. 
This implies $v^*u^ns^lv=0$ for all $n\in\Z$ and $l\in\{1,2,\ldots,k\}$, 
and hence completes the proof. 
Take $l\in\{1,2,\ldots,k\}$. 
We have 
$s^lg(u)=g(u^{K^{l}})s^l=g_l(u)s^l$ 
where $g_l\in C(\T)$ is defined 
by $g_l(z)=g(z^{K^{l}})$ for $z\in\T$. 
If $z\in\T$ satisfies $g(z)\neq 0$ 
then we get $z\in I_0$. 
Hence $z^{K^{l}}\in I_l$. 
Since $I_0\cap I_l=\emptyset$, 
we have $z^{K^{l}}\notin I_0$. 
Therefore $g_l(z)=g(z^{K^{l}})=0$. 
This shows $gg_l=0$. 
Thus we get $g(u)s^lg(u)=g(u)g_l(u)s^l=0$. 
We are done. 
\end{proof}

\begin{proposition}\label{Prop:axa=p}
Let $A,B\in M_N(\Z)$ satisfy the conditions 
{\rm (0)}, {\rm (1)} and {\rm (2)} in Section~\ref{Sec:OAB}. 
Then for every non-zero positive element $x\in \cO_{A,B}$, 
there exist $a\in \cO_{A,B}$ and 
$i\in \{1,2,\ldots,N\}$ 
with $a^*xa=p_i$. 
\end{proposition}

\begin{proof}
Take a non-zero positive element $x\in \cO_{A,B}$. 
To find $a\in \cO_{A,B}$ and 
$i\in \{1,2,\ldots,N\}$ 
with $a^*xa=p_i$, 
it suffices to find $a_0\in \cO_{A,B}$, 
$i\in \{1,2,\ldots,N\}$ and $C>0$ 
with $\|a_0^*xa_0-Cp_i\|<C$ 
because $p_i$ is a projection. 
Let $\e:=\|\varPhi(x)\|/3>0$ where 
$\varPhi\colon \cO_{A,B}\to \cO_{A,B}^{\T}$ 
is the faithful conditional expectation 
defined in Definition~\ref{Def:condexp}. 
Choose a positive element $x_0\in \acO_{A,B}$ with $\|x-x_0\|<\e$. 
Set $C:=\|\varPhi(x_0)\|$. 
We have 
\[
C=\|\varPhi(x_0)\|>\|\varPhi(x)\|-\e=2\e.
\]
Take $k\in\N$ with 
$x_0\in 
\spa\big\{s_{\xi}s_{\eta}^*\ \big|\ 
\xi,\eta\in\coprod_{l=0}^k \F^l\big\}$. 
Then $\varPhi(x_0)\in \Cr_{A,B}^k$. 
By Lemma~\ref{Lem:bab}, 
there exists $b\in \cspa\{s_{\xi}\mid \xi\in \F^k\}$ 
such that $\|b\|=1$, $\|b^*\varPhi(x_0)b\|=\|\varPhi(x_0)\|=C$ 
and $b^*\varPhi(x_0)b\in C^*(u_i)$ for some $i\in\{1,2,\ldots,N\}$. 
Take $f_0\in C(\T)$ with 
$f_0(u_i)=b^*\varPhi(x_0)b$ and $\|f_0\|=C$. 
Choose $f\in C(\T)$ 
which is $1$ on some non-empty open subset of $\T$, 
and satisfies $\|f_0-Cf\|<\e$. 
Set $s:=\s{0}{i,i}$. 
Since $B_{i,i}=1$, 
two elements $u_i$ and $s$ satisfy the assumptions 
of Lemma~\ref{Lem:isometry} for $K=A_{i,i}\geq 2$ 
in the \Ca generated by them. 
Applying Lemma~\ref{Lem:isometry} 
to $k\in\N$ and $f\in C(\T)$, 
we get $k_0\in\N$ with $k_0\geq k$ and $g\in C(\T)$ 
such that 
$v:=g(u_i)s^{k_0}\in \cO_{A,B}$ 
satisfies $v^*v=p_i$, $f(u_i)v=v$ 
and $v^*u_i^ns^lv=0$ for all $n\in\Z$ and $l\in\{1,2,\ldots,k\}$. 
We set $a_0:=b v$. 
Then we have $\|a_0\|\leq 1$. 
We will show $a_0^*s_{\xi}s_{\eta}^*a_0=0$ 
for $\xi\in\F^{l_1}$ and $\eta\in \F^{l_2}$ 
with $l_1,l_2\leq k$ and $l_1\neq l_2$. 
We may assume $l_1>l_2$. 
We set $l:=l_1-l_2\in \{1,2,\ldots,k\}$. 
By Lemma~\ref{Lem:sxsx*}, 
we see 
\[
b^*s_{\xi}s_{\eta}^*b\in \cspa\{s_{\zeta}\mid \zeta\in \F^l\}. 
\]
Hence it suffices to show $v^*s_{\zeta}v=0$ for all $\zeta\in \F^l$. 
Take $\mu\in \Omega_{A}^l$ with $\zeta\in F_\mu$. 
If $\mu\neq (i,i,\ldots,i)$ 
then $v^*s_{\zeta}=(s^{k_0})^*g(u_i)^*s_{\zeta}=0$ 
by Lemma~\ref{Lem:sxsx*}. 
If $\mu=(i,i,\ldots,i)$, 
then $s_{\zeta}=u_i^ns^l$ for some $n\in\Z$ 
because $B_{i,i}=1$. 
Hence we have $v^*s_{\zeta}v=0$. 
Thus we have shown $a_0^*s_{\xi}s_{\eta}^*a_0=0$ 
for $\xi\in\F^{l_1}$ and $\eta\in \F^{l_2}$ 
with $l_1,l_2\leq k$ and $l_1\neq l_2$. 
Therefore we obtain 
\[
a_0^*x_0a_0
=a_0^*\varPhi(x_0)a_0
=v^*f_0(u_i)v. 
\]
On the other hand, 
we have $v^*f(u_i)v=v^*v=p_i$. 
Thus 
we get $\|a_0^*x_0a_0-Cp_i\|<\e$. 
because $\|f_0-Cf\|<\e$. 
Hence 
\[
\|a_0^*xa_0-Cp_i\|
\leq \|a_0^*(x-x_0)a_0\|+\|a_0^*x_0a_0-Cp_i\|
<2\e<C.
\]
We are done. 
\end{proof}

\begin{proof}[Proof of Proposition~\ref{Prop:OAB_Ki}]
One can easily check that 
if $A,B$ satisfy the conditions {\rm (0)}, {\rm (1)} and {\rm (2)}, 
the projection $p_i\in \cO_{A,B}$ is full and 
properly infinite for all $i$. 
This fact and Proposition~\ref{Prop:axa=p} 
show that the \Ca $\cO_{A,B}$ 
is simple and purely infinite. 
\end{proof}

\section{$K$-theory of $\cO_{A,B}$}\label{Sec:K}

In this section, 
we give the proof of Proposition~\ref{Prop:KOAB} 
using the idea in \cite[Section~4]{Pi}. 
Let us take $N\in\fC$ and 
$A,B\in M_N(\Z)$ satisfying the condition (0). 

\begin{definition}\label{Def:TAB}
We define a \Ca $\cT_{A,B}$ to be the universal \Ca 
generated by mutually orthogonal projections $\{p_i\}_{i=1}^N$, 
partial unitaries $\{u_i\}_{i=1}^N$ with $u_i^0=p_i$, 
and partial isometries $\{\s{n}{i,j}\}_{(i,j)\in \Omega_{A},n\in\Z}$ 
satisfying the relations (i), (ii) 
in Definition~\ref{Def:OAB} and the relation 
\begin{itemize}
\item[(iii)'] 
$p_i\geq 
\sum_{j\in \Omega_{A}(i)}\sum_{n=1}^{A_{i,j}}\s{n}{i,j}\s{n}{i,j}^*$ 
for all $i$. 
\end{itemize}
\end{definition}

By the universality of $\cT_{A,B}$, 
there exists a natural surjection $\cT_{A,B}\to\cO_{A,B}$ 
whose kernel is denoted by $\cJ_{A,B}\subset \cT_{A,B}$. 
For each $i$, 
let us define $p_i',u_i'\in \cT_{A,B}$ by 
\begin{align*}
p_i'&:=
p_i-\sum_{j\in \Omega_{A}(i)}
\sum_{n=1}^{A_{i,j}}\s{n}{i,j}\s{n}{i,j}^*, \\
u_i'&:=
u_i-\sum_{j\in \Omega_{A}(i)}\sum_{n=1}^{A_{i,j}}
\s{n+B_{i,j}}{i,j}\s{n}{i,j}^*=u_ip_i'=p_i'u_i. 
\end{align*}
Then $p_i'$ is a projection and $u_i'$ is a partial unitary 
with $(u_i')^0=p_i'$. 
Note that $p_i',u_i'\in \cJ_{A,B}$ 
and $\cJ_{A,B}$ is generated by $\{p_i'\}_{i=1}^N$ as an ideal. 

\begin{lemma}\label{Lem:pisxi}
For $i\in\{1,2,\ldots,N\}$ and $\xi\in\F^*$, 
we have 
\[
p_i' s_\xi =
\begin{cases}
(u_i')^n & \text{if $\xi=n\xi_{(i)}^{0}\in F_{(i)}$ for $n\in\Z$,}\\
0 & \text{if $\xi\in F_{\mu}$ with $\mu\in \Omega_A^*\setminus\{(i)\}$.}
\end{cases}
\]
\end{lemma}

\begin{proof}
Straightforward. 
\end{proof}

\begin{lemma}\label{Lem:repTAB}
For each $i\in\{1,2,\ldots,N\}$, 
the spectrum of $u_i'\in \cJ_{A,B}\subset \cT_{A,B}$ contains $\T$. 
\end{lemma}

\begin{proof}
The proof is very similar to the one of Proposition~\ref{Prop:repOAB}. 
It suffices to find operators 
$\{P_i\}_{i=1}^N$, $\{U_i\}_{i=1}^N$, 
and $\{\bS{n}{i,j}\}_{(i,j)\in \Omega_{A},n\in\Z}$ 
on some Hilbert space 
such that they satisfy the relations (i), (ii) and (iii)' 
and the spectrum of 
\[
U_i':=U_i-
\sum_{j\in \Omega_{A}(i)}\sum_{n=1}^{A_{i,j}}
S(n+B_{i,j})_{i,j}\bS{n}{i,j}^* 
\]
contains $\T$ for each $i\in\{1,2,\ldots,N\}$. 

For $(i,j)\in \Omega_{A}$ 
and $\mu\in \Omega_{A}^k$ with $\mu_0=j$, 
we set $i\mu\in \Omega_{A}^{k+1}$ by $(i\mu)_0=i$ 
and $(i\mu)_l=\mu_{l-1}$ for $l=1,2,\ldots,k+1$. 
When $k=0$ i.e.\ $\mu=(j)$ and $i\mu=(i,j)$, 
we define a homomorphism 
$\sigma_{i}\colon F_\mu\to F_{i\mu}$ 
by $\sigma_{i}(\xi_{\mu}^{0})=A_{i,j}\xi_{i\mu}^{1}$. 
When $k$ is positive, 
we define a homomorphism $\sigma_{i}\colon F_\mu\to F_{i\mu}$ 
by $\sigma_{i}(\xi_{\mu}^{l})=\xi_{i\mu}^{l+1}$ for $l=1,2,\ldots,k$. 
By setting $\xi_{\mu}^{0}:=B_{\mu_0,\mu_1}\xi_{\mu}^{1}\in F_\mu$ 
for $\mu\in \Omega_{A}^k$ with positive $k$, 
we have $\sigma_{i}(\xi_{\mu}^{0})=A_{i,j}\xi_{\mu}^{1}$. 
Similarly as in the proof of Proposition~\ref{Prop:repOAB}, 
we can see 
that $\sigma_{i}\colon F_\mu\to F_{i\mu}$ 
is a well-defined injective homomorphism, 
and 
$F_{i\mu}
=\coprod_{n=1}^{A_{i,j}}\big(\sigma_{i}(F_\mu)+n\xi_{i\mu}^{1}\big)$.

We define operators $\{P_i\}_{i=1}^N$, $\{U_i\}_{i=1}^N$, 
and $\{\bS{n}{i,j}\}_{(i,j)\in \Omega_{A},n\in\Z}$ 
on the Hilbert space $\ell^2(\F^*)$ by 
\[
P_i(\delta_{\xi}):=
\begin{cases}
\delta_{\xi} & \text{if $\mu_0=i$}\\
0 & \text{if $\mu_0\neq i$},
\end{cases}\qquad 
U_i(\delta_{\xi}):=
\begin{cases}
\delta_{\xi+\xi_{\mu}^{0}} & \text{if $\mu_0=i$}\\
0 & \text{if $\mu_0\neq i$},
\end{cases}
\]
\[
\bS{n}{i,j}(\delta_{\xi}):=
\begin{cases}
\delta_{\sigma_{i}(\xi)+n\xi_{i\mu}^{1}} & \text{if $\mu_0=j$}\\
0 & \text{if $\mu_0\neq j$},
\end{cases}
\]
for $\xi\in F_\mu\subset\F^*$. 
Then it is routine to check that these operators satisfy 
the relations (i), (ii) and (iii)'. 
We see that for each $i\in\{1,2,\ldots,N\}$ 
the operator $U_i'\in B\big(\ell^2(\F^*)\big)$ 
acts as a bilateral shift on $\ell^2(F_{(i)})\cong\ell^2(\Z)$ 
and vanishes on $\ell^2(\F^*\setminus F_{(i)})$. 
Hence the spectrum of $U_i'\in B\big(\ell^2(\F^*)\big)$ contains $\T$. 
We are done. 
\end{proof}

\begin{remark}
By Corollary~\ref{Cor:injrepTAB}, 
we see that the representation of $\cT_{A,B}$ 
defined from the operators in the proof above 
is faithful. 
We use this fact in Appendix~\ref{Sec:cA}. 
\end{remark}

\begin{definition}
For $N\in\fC$, 
we define a \Ca $\cA_N$ to be the universal \Ca 
generated by mutually orthogonal projections $\{p_i\}_{i=1}^N$ 
and partial unitaries $\{u_i\}_{i=1}^N$ with $u_i^0=p_i$. 
\end{definition}

Clearly we have $\cA_N\cong C_0(\{1,2,\ldots,N\}\times \T)$. 
By definition, 
there exists a \shom $\cA_N\to \cT_{A,B}$ 
which sends $\{p_i,u_i\}_{i=1}^N$ to $\{p_i,u_i\}_{i=1}^N$. 
By Lemma~\ref{Lem:repTAB}, 
this \shom is injective. 
Hence we can consider $\cA_N$ as a \Csa of $\cT_{A,B}$. 
Note that the natural \shom $\cA_N\to \cO_{A,B}$ is injective 
only when $B$ has no zero rows by Corollary~\ref{Cor:spec(u_i)}. 

\begin{lemma}\label{Lem:cA}
The inclusion $\cA_N\hookrightarrow \cT_{A,B}$ 
is a $KK$-equivalence. 
\end{lemma}

\begin{proof}
See Appendix~\ref{Sec:cA}. 
\end{proof}

\begin{lemma}\label{Lem:cA'}
The \shom $\pi\colon \cA_N\to \cJ_{A,B}$ 
defined by $\pi(p_i)=p_i'$ and $\pi(u_i)=u_i'$ 
for $i=1,2,\ldots,N$ 
is an injection onto a hereditary and full \CsA . 
\end{lemma}

\begin{proof}
By Lemma~\ref{Lem:repTAB}, 
the \shom $\pi$ is an injection onto the \Csa of $\cJ_{A,B}$ 
generated by $\{p_i',u_i'\}_{i=1}^N$. 
This \Csa 
is the hereditary \Csa of $\cT_{A,B}$ 
generated by $\{p_i'\}_{i=1}^N$ by Lemma~\ref{Lem:pisxi}. 
Thus it is hereditary and full in $\cJ_{A,B}$ 
because $\cJ_{A,B}$ is generated by $\{p_i'\}_{i=1}^N$ as an ideal. 
\end{proof}

\begin{proposition}\label{Prop:OAB_UCT}
The \Ca $\cO_{A,B}$ is in the UCT class. 
\end{proposition}

\begin{proof}
Since the \Ca $\cA_N$ is commutative, 
Lemma~\ref{Lem:cA} and Lemma~\ref{Lem:cA'} show 
that $\cT_{A,B}$ and $\cJ_{A,B}$ are in the UCT class 
(see \cite[Definition~2.4.5]{RS}). 
Now the ``two out of three principle'' 
(\cite[Proposition~2.4.7~(i)]{RS}) shows 
that $\cO_{A,B}\cong \cT_{A,B}/\cJ_{A,B}$ is in the UCT class. 
\end{proof}

\begin{lemma}\label{Lem:I->T}
The homomorphisms $K_0(\cJ_{A,B})\to K_0(\cT_{A,B})$ and 
$K_1(\cJ_{A,B})\to K_1(\cT_{A,B})$ 
induced by the embedding $\cJ_{A,B}\hookrightarrow \cT_{A,B}$ 
coincide with $I-A\colon \Z^N\to \Z^N$ and 
$I-B\colon \Z^N\to \Z^N$, respectively. 
\end{lemma}

\begin{proof}
By Lemma~\ref{Lem:cA} and Lemma~\ref{Lem:cA'}, 
all of the four abelian groups $K_0(\cT_{A,B})$, $K_1(\cT_{A,B})$, 
$K_0(\cJ_{A,B})$ and $K_1(\cJ_{A,B})$ are isomorphic to $\Z^N$ 
whose basis are given by $\{[p_i]\}_{i=1}^N$, $\{[u_i]\}_{i=1}^N$, 
$\{[p_i']\}_{i=1}^N$ and $\{[u_i']\}_{i=1}^N$, respectively. 

For each $i$, 
we have 
\[
[p_i']= [p_i] - \sum_{j\in \Omega_{A}(i)}
\sum_{n=1}^{A_{i,j}}[\s{n}{i,j}\s{n}{i,j}^*]
=[p_i] - \sum_{j=1}^N A_{i,j}[p_j]. 
\]
This shows the result on $K_0$. 
For $(i,j)\in\Omega_A$, we define a partial unitary 
\[
u_{i,j}:=\sum_{n=1}^{A_{i,j}}\s{n+1}{i,j}\s{n}{i,j}^*\in \cO_{A,B}.
\]
Then we have $[u_{i,j}]=[u_j]$ (see Lemma~\ref{Lem:DefU}). 
Hence for each $i$, 
we get 
\[
[u_i']
=[u_i]-\sum_{j\in \Omega_{A}(i)} [u_{i,j}^{B_{i,j}}]
=[u_i]-\sum_{j=1}^N B_{i,j}[u_j] 
\]
(see Subsection~\ref{ssec:partunit}). 
Thus we get the conclusion for $K_1$. 
\end{proof}

\begin{proposition}\label{Prop:exactseqs}
We have $\Gamma_{A,B}$-equivariant exact sequences 
\[
\begin{CD}
0 @>>> \coker (I-A) @>>> K_0(\cO_{A,B}) @>>> \ker (I-B) @>>> 0,\\
0 @>>> \coker (I-B) @>>> K_1(\cO_{A,B}) @>>> \ker (I-A) @>>> 0\phantom{,}\\
\end{CD}
\]
both of which have $\Gamma_{A,B}$-equivariant splitting maps. 
\end{proposition}

\begin{proof}
The same formula as $\Gamma_{A,B}\act \cO_{A,B}$ 
defines an action $\Gamma_{A,B}\act \cT_{A,B}$ 
which makes the short exact sequence 
\[
\begin{CD}
0 @>>> \cJ_{A,B} @>>> \cT_{A,B} @>>> \cO_{A,B} @>>> 0
\end{CD}
\]
$\Gamma_{A,B}$-equivariant. 
From this sequence, 
we get a $\Gamma_{A,B}$-equivariant $6$-term exact sequence 
\[
\begin{CD}
K_0(\cJ_{A,B}) @>>> K_0(\cT_{A,B}) @>>> K_0(\cO_{A,B}) \\
@AAA @. @VVV \\
K_1(\cO_{A,B}) @>>> K_1(\cT_{A,B}) @>>> K_1(\cJ_{A,B}). 
\end{CD}
\]
By Lemma~\ref{Lem:I->T} and its proof, 
the homomorphisms $K_0(\cJ_{A,B})\to K_0(\cT_{A,B})$ and 
$K_1(\cJ_{A,B})\to K_1(\cT_{A,B})$ 
coincide with $I-A\colon \Z^N\to \Z^N$ and 
$I-B\colon \Z^N\to \Z^N$, $\Gamma_{A,B}$-equivariantly. 
Thus we get the desired $\Gamma_{A,B}$-equivariant short exact sequences. 
Since $\ker (I-A)$ and $\ker (I-B)$ are free abelian groups, 
these sequences split. 
That we can find $\Gamma_{A,B}$-equivariant splitting maps 
will be proven in Appendix~\ref{Sec:split}. 
\end{proof}

\begin{remark}\label{Rem:Ker=0}
In the proof of main theorem in Section~\ref{Sec:MainThm}, 
we only use Proposition~\ref{Prop:exactseqs} 
for the case $\ker (I-A)=\ker (I-B)=0$ 
and hence we do not need the results on 
$\Gamma_{A,B}$-equivariant splitting maps. 
However in order to compute examples, 
it is important to show that there exist 
$\Gamma_{A,B}$-equivariant splitting maps. 
\end{remark}

\appendix

\section{Proof of Lemma~\ref{Lem:cA}}\label{Sec:cA}

In this appendix, 
we give a proof of Lemma~\ref{Lem:cA} in an explicit way, 
following the argument in \cite[Appendix~C]{Ka3} 
which was inspired by the original proof in \cite{Pi}. 
We first need the following proposition 
which has own importance. 

\begin{proposition}\label{Prop:essential}
For $A,B\in M_N(\Z)$ satisfying the condition {\rm (0)}, 
the ideal $\cJ_{A,B}$ of $\cT_{A,B}$ is essential. 
\end{proposition}

\begin{proof}
We use the notation established in the proof of Lemma~\ref{Lem:Ck}. 

Let $\cJ_{A,B}^\perp$ be the ideal of $\cT_{A,B}$ defined by 
\[
\cJ_{A,B}^\perp
:=\{x\in \cT_{A,B}\mid \text{$xy=0$ for all $y\in \cJ_{A,B}$}\}.
\]
We will show $\cJ_{A,B}^\perp=0$. 
The same formula of the gauge action 
$\beta\colon \T\act \cO_{A,B}$ defined 
in Definition~\ref{Def:gauge} 
defines an action $\T\act \cT_{A,B}$. 
Since $\cJ_{A,B}$ is invariant under this action, 
so is $\cJ_{A,B}^\perp$. 
Hence in order to show $\cJ_{A,B}^\perp=0$, 
it suffices to see $\cJ_{A,B}^\perp\cap \cT_{A,B}^\T=0$ 
where $\cT_{A,B}^\T$ is 
the fixed point algebra of the action $\T\act \cT_{A,B}$. 

For each $k\in\N$, 
we define \CsA s $\Cr_{A,B}^k,\D_{A,B}^k$ of $\cT_{A,B}$ by 
\begin{align*}
\Cr_{A,B}^k
:=\cspa{}&\big\{s_{\xi}s_{\eta}^*\ \big|\ \xi,\eta\in \F^k\big\} \\
=\cspa{}&\big\{s_{\xi}p_is_{\eta}^*\ \big|\ \xi,\eta\in \F^k_i, 
i\in\{1,2,\ldots,N\}\big\}, \\
\D_{A,B}^k
:=\cspa{}&\big\{s_{\xi}p_i's_{\eta}^*\ \big|\ \xi,\eta\in \F^k_i, 
i\in\{1,2,\ldots,N\}\big\}. 
\end{align*}
Similarly as in Lemma~\ref{Lem:indlim}, 
we have 
\[
\cT_{A,B}^{\T}=\overline{\bigcup_{k\in\N}\Big(\sum_{l=0}^k \Cr_{A,B}^l\Big)}. 
\]
Hence in order to show $\cJ_{A,B}^\perp\cap \cT_{A,B}^\T=0$, 
it suffices to see $\cJ_{A,B}^\perp \cap (\sum_{l=0}^k \Cr_{A,B}^l)=0$ 
for all $k\in\N$. 
This reduces to the problem showing $\cJ_{A,B}^\perp \cap \Cr_{A,B}^k=0$ 
because the \Ca $\sum_{l=0}^k \Cr_{A,B}^l$ has the orthogonal decomposition 
\[
\sum_{l=0}^k \Cr_{A,B}^l
=\sum_{l=0}^{k-1} \D_{A,B}^l + \Cr_{A,B}^k. 
\]
and $\D_{A,B}^l\subset \cJ_{A,B}$ for all $l$. 

Fix $k\in\N$, 
and we will show $\cJ_{A,B}^\perp \cap \Cr_{A,B}^k=0$. 
Similarly as in the proof of Lemma~\ref{Lem:Ck}, 
we obtain 
\begin{align*}
\Cr_{A,B}^k 
&\cong \bigoplus_{i=1}^N C^*(u_i)\otimes \cK\big(\ell^2(X^k_i)\big),&
\D_{A,B}^k 
&\cong \bigoplus_{i=1}^N C^*(u_i')\otimes \cK\big(\ell^2(X^k_i)\big).
\end{align*}
By Lemma~\ref{Lem:repTAB}, 
$C^*(u_i)\cong C^*(u_i')\cong C(\T)$ for each $i$. 
This observation and the proof of Lemma~\ref{Lem:Ck} 
imply that 
the \shom $\varphi\colon \Cr_{A,B}^k\to \D_{A,B}^k$ defined by 
\[
\varphi(s_{\xi}p_is_{\eta}^*):=s_{\xi}p_i's_{\eta}^*
\]
for $\xi,\eta\in \F^k_i$ and $i\in\{1,2,\ldots,N\}$ 
is an isomorphism. 
For each $i\in\{1,2,\ldots,N\}$, 
we define a projection $p_i''\in\cT_{A,B}$ by 
\[
p_i''
:=p_i-\sum_{\mu\in \Omega_A^{k+1}(i)}\ \sum_{\xi\in X_\mu}s_\xi s_\xi^*
\]
where $\Omega_A^{k+1}(i)$ is defined by 
\[
\Omega_A^{k+1}(i)
:=\{\mu\in \Omega_A^{k+1}\mid \mu_0=i\} 
\]
which is a finite set by the condition (0). 
We have $p_i''\in\cJ_{A,B}$ 
because the image of $p_i''\in \cT_{A,B}$ in $\cO_{A,B}$ 
can be shown to be $0$. 
We have $p_{\mu_0}''s_\xi=s_\xi p_{\mu_k}'$ for 
$\xi\in F_\mu\subset \F^k$. 
This implies that 
the isomorphism $\varphi$ satisfies 
\[
\varphi(x)=\sum_{i=1}^N p_{i}''x=\sum_{i=1}^N xp_{i}''
\]
for $x\in \Cr_{A,B}^k$. 
Since $p_i''\in\cJ_{A,B}$ for all $i$, 
$\Cr_{A,B}^k\cap \cJ_{A,B}^\perp$ 
is contained in the kernel of $\varphi$. 
Thus we have $\Cr_{A,B}^k\cap \cJ_{A,B}^\perp=0$. 
This completes the proof. 
\end{proof}

\begin{corollary}\label{Cor:injrepTAB}
A \shom $\varphi$ from the \Ca $\cT_{A,B}$ is injective 
if and only if 
the spectrum of $\varphi(u_i')$ contains $\T$ 
for every $i\in \{1,2,\ldots,N\}$. 
\end{corollary}

\begin{proof}
By Proposition~\ref{Prop:essential}, 
$\varphi$ is injective 
if and only if its restriction to $\cJ_{A,B}$ is injective. 
This happens exactly when 
the restriction of $\varphi$ to the image of 
the \shom $\pi\colon \cA_N\to \cJ_{A,B}$ in Lemma~\ref{Lem:cA'} is injective. 
This condition is equivalent that 
the spectrum of $\varphi(u_i')$ contains $\T$ 
for every $i\in \{1,2,\ldots,N\}$. 
\end{proof}

Recall that for two (separable) \CA s $\cA$ and $\cB$, 
the $KK$-group $KK(\cA,\cB)$ is an abelian group such that 
a homotopy equivalence class of a \shom $\varphi\colon \cA \to \cB$ 
defines an element $[\varphi]$ of $KK(\cA,\cB)$ 
(see \cite[Chapter~VIII]{Bl} for definitions and results stated below). 
We say that a \shom $\varphi\colon \cA \to \cB$ is 
a {\em $KK$-equivalence} 
if the element $[\varphi]\in KK(\cA,\cB)$ 
has an inverse in $KK(\cB,\cA)$ with respect to 
the Kasparov product. 
The facts on $KK$-equivalences we use in the following argument 
are summarized as follows: 
\begin{itemize}
\item A composition of two $KK$-equivalences is a $KK$-equivalence. 
\item Let $\varphi_i\colon \cA_{i}\to \cA_{i+1}$ for $i=1,2,3$ 
be $*$-ho\-mo\-mor\-phisms. 
If both $\varphi_2\circ \varphi_1$ and $\varphi_3\circ \varphi_2$ 
are $KK$-equivalences, 
then all the three \shoms $\varphi_1, \varphi_2, \varphi_3$ 
are $KK$-equivalences. 
\item An inclusion map of a hereditary and full \Csa is a $KK$-equivalence. 
\item Let $0\to \cI\to \cA\to \cB\to 0$ be a splitting short exact sequence. 
If $\cB$ is contractible, 
then the map $\cI\to \cA$ is a $KK$-equivalence. 
\item A \shom $\cA\to \cB$ is a $KK$-equivalence if and only if 
the induced \shom $S\cA\to S\cB$ is a $KK$-equivalence. 
\end{itemize}
Here we define $S\cA=C_0\big((0,1),\cA\big)$ 
for a \Ca $\cA$. 
A $KK$-equivalence induces isomorphisms of $K$-groups, 
and conversely a \shom $\cA\to \cB$ 
inducing isomorphisms of $K$-groups is a $KK$-equivalence 
if $\cA$ and $\cB$ are in the UCT class. 
It is easy to see that the five statements above are still valid 
if we replace ``$KK$-equivalence'' to 
``\shom inducing isomorphisms of $K$-groups''. 

Let us take $A,B\in M_N(\Z)$ satisfying the condition (0). 
We will show that the inclusion $\cA_N\hookrightarrow \cT_{A,B}$ 
is a $KK$-equivalence. 
For $k,l\in\N$, 
we define 
\[
\Omega_A^{k,l}
:=\big\{(\mu,\mu')\ \big|\ \text{$\mu\in \Omega_A^k$, 
$\mu'\in \Omega_A^l$ with $\mu_k=\mu'_0$} \big\}. 
\]
Note that the map 
\[
\Omega_A^{k,l}\ni (\mu,\mu')\mapsto \mu\mu'\in \Omega_A^{k+l}
\]
is bijective where $\mu\mu'\in \Omega_A^{k+l}$ is defined by 
$(\mu\mu')_m=\mu_m$ for $m=0,1,\ldots, k$ 
and $(\mu\mu')_{k+m}=\mu'_m$ for $m=0,1,\ldots, l$. 
For $(\mu,\mu')\in \Omega_A^{k,l}$, 
we set $F_{\mu,\mu'}:=F_{\mu\mu'}$. 
Note that there exists a homomorphism $F_\mu\to F_{\mu,\mu'}$ 
which sends $\xi_\mu^k$ to $\xi_{\mu\mu'}^k$. 
We set $\Omega_A^{*,*}:=\coprod_{k,l\in\N}\Omega_A^{k,l}$ and  
$\F^{*,*}:=\coprod_{(\mu,\mu')\in \Omega_A^{*,*}}F_{\mu,\mu'}$. 

Take $(i,j)\in \Omega_{A}$ 
and $(\mu,\mu')\in \Omega_A^{k,l}$ with $\mu_0=j$. 
Note that the element $i\mu\in \Omega_A^{k+1}$ 
defined in the proof of Lemma~\ref{Lem:repTAB} 
is nothing but $(i,j)\mu$ defined here, 
and that we have $(i\mu)\mu'=i(\mu\mu')$. 
We denote by the same symbol $\sigma_{i}$ 
the homomorphism $F_{\mu,\mu'}\to F_{i\mu,\mu'}$ 
induced by the homomorphism 
$\sigma_{i}\colon F_{\mu\mu'}\to F_{i\mu\mu'}$ 
defined in the proof of Lemma~\ref{Lem:repTAB}. 
We define a \shom $\rho\colon\cT_{A,B}\to B\big(\ell^2(\F^{*,*})\big)$ by 
\begin{align*}
\rho(p_i)(\delta_{\xi})&:=
\begin{cases}
\delta_{\xi} & \text{if $\mu_0=i$}\\
0 & \text{if $\mu_0\neq i$},
\end{cases}\qquad
\rho(u_i)(\delta_{\xi}):=
\begin{cases}
\delta_{\xi+\xi_{\mu\mu'}^{0}} & \text{if $\mu_0=i$}\\
0 & \text{if $\mu_0\neq i$},
\end{cases}\\
\rho(\s{n}{i,j})(\delta_{\xi})&:=
\begin{cases}
\delta_{\sigma_{i}(\xi)+n\xi_{i\mu\mu'}^{1}} & \text{if $\mu_0=j$}\\
0 & \text{if $\mu_0\neq j$},
\end{cases}
\end{align*}
for generators $\{p_i,u_i,\s{n}{i,j}\}$ of $\cT_{A,B}$ and 
$\xi\in F_{\mu,\mu'}\subset \F^{*,*}$. 
It is routine to check that this is well-defined. 

Since we can naturally identify $\Omega_A^{0,l}$ and $\Omega_A^{l}$, 
we consider $\Omega_A^{*}$ as a subset of $\Omega_A^{*,*}$ 
by $\Omega_A^{*}=\coprod_{l\in\N}\Omega_A^{0,l}\subset \Omega_A^{*,*}$. 
Using this identification, 
we can consider $\F^*$ as a subset of $\F^{*,*}$. 
Then the representation 
defined from the operators in the proof of Lemma~\ref{Lem:repTAB} 
induces the representation 
$\cT_{A,B}\to B\big(\ell^2(\F^{*})\big) \subset B\big(\ell^2(\F^{*,*})\big)$ 
which is denoted by $\tpi$. 
By Corollary~\ref{Cor:injrepTAB}, 
the representation $\tpi$ is faithful. 
It is routine to check $\tpi(p_i)=\rho(p_i')$ 
and $\tpi(u_i)=\rho(u_i')$ for $i\in\{1,2,\ldots,N\}$. 
Let $\cL\subset B\big(\ell^2(\F^{*,*})\big)$ 
be the \Ca generated by $\tpi(\cT_{A,B})\cup\rho(\cT_{A,B})$. 
For $\xi\in\F^k$ with $k\geq 1$ and $i\in \{1,2,\ldots,N\}$, 
we have $\tpi(p_i)\rho(s_\xi)=\rho(p'_i s_{\xi})=0$ 
by Lemma~\ref{Lem:pisxi}. 
This fact shows 
that $\tpi(\cT_{A,B})\subset \cL$ 
is the hereditary \Csa of $\cL$ 
generated by $\{\tpi(p_i)\}_{i=1}^N$. 
We let $\cK\subset \cL$ be the ideal of $\cL$ 
generated by $\{\tpi(p_i)\}_{i=1}^N$. 
We have $\rho(\cJ_{A,B})\subset \cK$ 
because $\rho(p_i')=\tpi(p_i)$ for all $i$. 

Let $\pi\colon\cA_N\to \cJ_{A,B}$ be the \shom 
defined in Lemma~\ref{Lem:cA'}. 
Then the restriction of $\tpi$ to $\cA_N\subset \cT_{A,B}$ 
coincides with $\rho\circ\pi$. 
Thus we get the following commutative diagram: 
\[
\SelectTips{cm}{}
\xymatrix{
\cA_N \ar@<-0.2ex>@{^{(}->}[d] \ar[r]^\pi 
& \cJ_{A,B} \ar[d]^\rho \ar@<-0.2ex>@{^{(}->}[r]
& \cT_{A,B} \ar[d]^\rho \\
\cT_{A,B} \ar[r]^{\tpi} 
& \cK \ar@<-0.2ex>@{^{(}->}[r] 
& \cL}
\]
We define 
\begin{align*}
\cD&:=\big\{f\in C_0\big((-1,1),\cT_{A,B}\big)\ \big|\ 
f(s)-f(-s)\in \cJ_{A,B}\mbox{ for all }s\in (-1,1)\big\},\\
\tcD&:=\big\{f\in C_0\big((-1,1),\cL\big)\ \big|\ 
f(s)-f(-s)\in \cK\mbox{ for all }s\in (-1,1)\big\}.
\end{align*}
Then $\rho$ induces the \shom $\cD\to \tcD$ which is also denoted by $\rho$. 
We have the following splitting short exact sequences; 
\[
\begin{CD}
0 @>>> S\cJ_{A,B} @>>> \cD @>>> C_0\big((-1,0],\cT_{A,B}\big) @>>> 0,\\
0 @>>> S\cK @>>> \tcD @>>> C_0\big((-1,0],\cL\big) @>>> 0.\\
\end{CD}
\]
Since $C_0\big((-1,0],\cT_{A,B}\big)$ is contractible, 
the inclusion map $S\cJ_{A,B}\hookrightarrow \cD$ 
is a $KK$-equivalence. 
The \shom $S\cA_N\to S\cJ_{A,B}$ induced by the \shom $\pi$ 
is a $KK$-equivalence
because so is $\pi$ by Lemma~\ref{Lem:cA'}. 
Hence the composition of these two \shoms 
which is denoted by $\varphi\colon S\cA_N\to \cD$ 
is a $KK$-equivalence. 
Since the map $\tpi\colon \cT_{A,B}\to \cK$ is 
an injection onto a full and hereditary \Csa of $\cK$, 
we can show in a very similar way as above 
that the \shom $\tphi\colon S\cT_{A,B}\to \tcD$ induced by $\tpi$ 
is a $KK$-equivalence. 
We get $\tphi\circ\iota=\rho\circ\varphi$ 
where $\iota\colon S\cA_N\hookrightarrow S\cT_{A,B}$ 
is the embedding map. 
\[
\SelectTips{cm}{}
\xymatrix{
S\cA_N \ar@<-0.3ex>@{^{(}->}[d]^\iota \ar[r]^\varphi 
& \cD \ar[d]^\rho \\
S\cT_{A,B} \ar@{.>}[ur]|{\psi} \ar[r]^{\tphi} 
& \tcD}
\]
We will construct a \shom $\psi\colon S\cT_{A,B}\to \cD$ 
and show that $\psi\circ\iota$ and $\rho\circ\psi$ 
induces the same elements as $\varphi$ and $\tphi$ 
in $KK$-groups, respectively. 
From this fact, 
we see that the map 
$\iota\colon S\cA_N\hookrightarrow S\cT_{A,B}$, 
and hence the inclusion $\cA_N\hookrightarrow \cT_{A,B}$ 
is a $KK$-equivalence. 

We define an endomorphism $\theta\colon \cT_{A,B}\to\cT_{A,B}$ by 
\begin{align*}
\theta(p_i)
&:=\sum_{j\in \Omega_{A}(i)}
\sum_{n=1}^{A_{i,j}}\s{n}{i,j}\s{n}{i,j}^*=p_i-p_i'\\
\theta(u_i)
&:=\sum_{n=1}^{A_{i,j}}\s{n+B_{i,j}}{i,j}\s{n}{i,j}^*
=u_i-u_i'=u_i\theta(p_i)\\
\theta(\s{n}{i,j})
&:=\s{n}{i,j}\theta(p_j)
\end{align*}
for generators $\{p_i,u_i,\s{n}{i,j}\}$ of $\cT_{A,B}$. 
We can show that it is well-defined and 
satisfies $x-\theta(x)\in \cJ_{A,B}$ for all $x\in \cT_{A,B}$ 
because this holds for the generators. 
Hence we can define a \shom $\psi\colon S\cT_{A,B}\to \cD$ by 
\[
\psi(f)(s)=
\begin{cases}
f(s) & \text{for $s\in [0,1)$,}\\
\theta\big(f(-s)\big) & \text{for $s\in (-1,0]$}
\end{cases}
\]
for $f\in S\cT_{A,B}$. 
The proof of Lemma~\ref{Lem:cA} completes once we see that 
$\psi\circ\iota$ and $\rho\circ\psi$ induce 
the same elements in $KK$-groups as 
$\varphi$ and $\tphi$, respectively.

Let $\varphi'\colon S\cA_N\to \cD$ be the \shom defined by 
\[
\varphi'(f)(s)=
\begin{cases}
\theta\big(f(s)\big) & \text{for $s\in [0,1)$,}\\
\theta\big(f(-s)\big) & \text{for $s\in (-1,0]$.}
\end{cases}
\]
Since $\pi'+\theta\colon \cA_N\to \cT_{A,B}$ coincides with the embedding map, 
we have $\psi\circ\iota=\varphi+\varphi'$. 
The element of $KK(S\cA_N,\cD)$ 
induced by the \shom $\varphi'$ is $0$ 
because it factors through 
the contractible \Ca 
$C_0\big([0,1),\cA_N\big)$. 
Hence the element of $KK(S\cA_N,\cD)$ induced by $\psi\circ\iota$ 
coincides with the one by $\varphi$. 

Since the representation $\rho\circ\theta$ vanishes 
on $\ell^2(\F^{*})$, 
$\tpi+(\rho\circ\theta)$ is a representation of $\cT_{A,B}$ 
on $\ell^2(\F^{*,*})$. 
The two representations $\tpi+(\rho\circ\theta)$ and $\rho$ coincide 
for $\{p_i,u_i\}_{i=1}^N$, 
but do not coincide 
for $\{\s{n}{i,j}\}_{(i,j)\in \Omega_{A},n\in\Z}$. 
These two representations can be connected via 
a continuous path 
$\{\rho_t\}_{t\in [0,1]}$ 
of representations defined by 
$\rho_t(p_i)=\rho(p_i)$, 
$\rho_t(u_i)=\rho(u_i)$ 
and 
\[
\rho_t(\s{n}{i,j}):=\sqrt{1-t^2}\tpi(\s{n}{i,j})
+ t\rho(\s{n}{i,j}p_j') + \rho(\theta(\s{n}{i,j})) 
\]
for generators $\{p_i,u_i,\s{n}{i,j}\}$ of $\cT_{A,B}$ and $t\in [0,1]$. 
It is routine to check that for each $t\in [0,1]$, 
the \shom $\rho_t\colon \cT_{A,B}\to\cL$ is 
well-defined and satisfies 
$\rho_t(x)-\rho(\theta(x))\in \cK$ 
for all $x\in \cT_{A,B}$. 
We define a homotopy 
$\{\tphi_t\}_{t\in [0,1]}$ 
of \shoms from $S\cT_{A,B}$ to $\tcD$ by 
\[
\tphi_t(f)(s)=
\begin{cases}
\rho_t(f(s)) & \text{for $s\in [0,1)$,}\\
\rho\big(\theta(f(-s))\big) & \text{for $s\in (-1,0]$}
\end{cases}
\]
Similarly as above, 
we can see that the two elements $\tphi$ and $\tphi_0$ 
induce the same element in $KK(S\cT_{A,B},\tcD)$ 
because $\tphi_0$ is the sum of $\tphi$ 
and a \shom factoring through a contractible \CA . 
Since $\tphi_0$ and $\tphi_1=\rho\circ\psi$ also 
induce the same element, 
the element in $KK(S\cT_{A,B},\tcD)$ 
induced by $\tphi$ coincides with the one induced by $\rho\circ\psi$. 
This completes the proof of Lemma~\ref{Lem:cA}.

\section{$\Gamma_{A,B}$-equivariant splitting maps}\label{Sec:split}

In this appendix, 
we prove that 
the $\Gamma_{A,B}$-equivariant exact sequences 
\[
\begin{CD}
0 @>>> \coker (I-A) @>>> K_0(\cO_{A,B}) @>>> \ker (I-B) @>>> 0,\\
0 @>>> \coker (I-B) @>>> K_1(\cO_{A,B}) @>>> \ker (I-A) @>>> 0\phantom{,}\\
\end{CD}
\]
in Proposition~\ref{Prop:exactseqs} 
have $\Gamma_{A,B}$-equivariant splitting maps. 

We first see what the homomorphism $K_0(\cO_{A,B})\to \ker (I-B)$ is. 
An element in $K_0(\cO_{A,B})$ can be presented 
by a unitary $u\in C\big([0,1],M_n(\widetilde{\cO}_{A,B})\big)$ 
with $u(0)=u(1)=1$ for some positive integer $n$, 
where $\widetilde{\cO}_{A,B}$ is 
the minimal unital \Ca containing $\cO_{A,B}$. 
We can lift $u\in C\big([0,1],M_n(\widetilde{\cO}_{A,B})\big)$ 
to a unitary $\bar{u}\in C\big([0,1],M_n(\widetilde{\cT}_{A,B})\big)$ 
with $\bar{u}(0)=1$. 
Then $\bar{u}(1)$ is a unitary in $M_n(\widetilde{\cJ}_{A,B})$ 
which defines an element in $K_1(\cJ_{A,B})$. 
In this way, 
we have the homomorphism $K_0(\cO_{A,B})\to K_1(\cJ_{A,B})\cong \Z^N$, 
and by the proof of Proposition~\ref{Prop:exactseqs} 
we can see that 
the image of this homomorphism is $\ker (I-B)\subset \Z^N$. 

For $(i,j)\in\Omega_A$, 
we set $u_{i,j}
=\sum_{n=1}^{A_{i,j}}\s{n+1}{i,j}\s{n}{i,j}^*\in \cT_{A,B}$. 
Note that 
$\{u_i'\}_{i=1}^N \cup \{u_{i,j}\}_{(i,j)\in\Omega_A}$ 
is a family of mutually orthogonal partial unitaries in $\cT_{A,B}$. 
Hence $\{\tu_i'\}_{i=1}^N \cup \{\tu_{i,j}\}_{(i,j)\in\Omega_A}$ 
is a family of mutually commutative unitaries in $\tcT_{A,B}$. 
We can see that these unitaries commute 
with mutually commutative unitaries $\{\tu_i\}_{i=1}^N$ 
from the equations 
$\tu_i=\tu_i'\times\prod_{j\in \Omega_A(i)}\tu_{i,j}^{B_{i,j}}$. 

\begin{lemma}\label{Lem:uu^n}
For $\sum_{i=1}^N n_ie_i\in \ker (I-B)$, 
we have 
$\prod_{(i,j)\in\Omega_A}\big(\tu_j \tu_{i,j}^{-1}\big)^{n_iB_{i,j}}
=\prod_{i=1}^N (\tu'_i)^{n_i}$. 
\end{lemma}

\begin{proof}
Take $\sum_{i=1}^N n_ie_i\in \ker (I-B)$. 
Then for each $j\in\{1,2,\ldots,N\}$, 
we have 
$n_j=\sum_{i=1}^Nn_iB_{i,j}$. 
Hence we obtain 
\begin{align*}
\prod_{(i,j)\in\Omega_A}\big(\tu_j \tu_{i,j}^{-1}\big)^{n_iB_{i,j}}
&= \prod_{i,j=1}^N \tu_j^{n_iB_{i,j}}\times 
\prod_{i=1}^N\Big(\prod_{j\in\Omega_A(i)}\tu_{i,j}^{B_{i,j}}\Big)^{-n_i}\\
&= \prod_{j=1}^N \tu_j^{\sum_{i=1}^Nn_iB_{i,j}}\times 
\prod_{i=1}^N\big(\tu_i(\tu_i')^{-1}\big)^{-n_i}\\
&= \prod_{j=1}^N \tu_j^{n_j}\times 
\prod_{i=1}^N(\tu_i^{-n_i}(\tu_i')^{n_i})\\
&=\prod_{i=1}^N (\tu'_i)^{n_i}. 
\end{align*}
\end{proof}

The action $\Gamma_{A,B}\act \cT_{A,B}$ naturally induces 
an  action $\Gamma_{A,B}\act C\big([0,1],M_2(\tcT_{A,B})\big)$. 

\begin{lemma}\label{Lem:DefU}
There exists a family of unitaries $\{U_{i,j}\}_{(i,j)\in \Omega_A}$ 
in $C\big([0,1],M_2(\tcT_{A,B})\big)$ such that 
\[
U_{i,j}(0)=\mat{1}{0}{0}{1},\qquad 
U_{i,j}(1)=\mat{\tu_j\tu_{i,j}^{-1}}{0}{0}{1}
\]
and $\gamma(U_{i,j})=U_{\gamma(i),\gamma(j)}$ 
for $(i,j)\in \Omega_A$ and $\gamma\in\Gamma_{A,B}$. 
\end{lemma}

\begin{proof}
Take $(i,j)\in \Omega_A$. 
We define a projection $q_{i,j}\in \cT_{A,B}$ 
and a partial unitary $v_{i,j}\in \cT_{A,B}$ 
with $v_{i,j}^0=q_{i,j}$ by 
\[
q_{i,j}
:=\sum_{n=1}^{A_{i,j}}\s{n}{i,j}\s{n}{i,j}^*,\qquad
v_{i,j}
:=\s{A_{i,j}}{i,j}\s{1}{i,j}^*
+\sum_{n=2}^{A_{i,j}}\s{n-1}{i,j}\s{n}{i,j}^*.
\]
We have $v_{i,j}^{A_{i,j}}=q_{i,j}$. 
For each $n=0,1,\ldots A_{i,j}-1$, 
the spectrum projection $q_{i,j,n}$ of $v_{i,j}$ 
corresponding to the spectrum $e^{2\pi\sqrt{-1}n/A_{i,j}}$ of $v_{i,j}$ 
is computed as 
\[
q_{i,j,n}=
\frac{1}{A_{i,j}}
\sum_{m=0}^{A_{i,j}-1}
\big(
e^{2\pi\sqrt{-1}\frac{-n}{A_{i,j}}}
v_{i,j}\big)^m. 
\]
Then the homotopy $[0,1]\ni t\mapsto v_{i,j}(t)\in \cT_{A,B}$ 
defined by 
\[
v_{i,j}(t):=
\sum_{n=0}^{A_{i,j}-1}
e^{2\pi\sqrt{-1}\frac{nt}{A_{i,j}}}
q_{i,j,n}
\]
satisfies $v_{i,j}(0)=q_{i,j}$, $v_{i,j}(1)=v_{i,j}$ 
and $v_{i,j}(t)$ is a partial unitary 
with $v_{i,j}(t)^0=q_{i,j}$. 
We define a homotopy 
$[0,1]\ni t\mapsto W_{i,j}(t)\in C\big([0,1],M_2(\cT_{A,B})\big)$ 
of partial unitaries by 
\[
W_{i,j}(t):=
V_{i,j}(t) \mat{u_{i,j}v_{i,j}(t)}{0}{0}{p_{j}} V_{i,j}(t)^* 
\]
where 
\[
V_{i,j}(t):=
\mat{\sqrt{1-t^2}\s{1}{i,j}(\s{1}{i,j})^*}
{-t\s{1}{i,j}}{t(\s{1}{i,j})^*}{\sqrt{1-t^2}p_j}
+\mat{\sum_{n=2}^{A_{i,j}}\s{n}{i,j}\s{n}{i,j}^*}{0}{0}{0} 
\]
is also a partial unitary. 
We can compute 
\[
W_{i,j}(0)=\mat{u_{i,j}}{0}{0}{p_j},\qquad 
W_{i,j}(1)=\mat{q_{i,j}}{0}{0}{u_j}
\]
using the equation 
\[
u_{i,j}v_{i,j}
=\s{1}{i,j}u_j\s{1}{i,j}^*
+\sum_{n=2}^{A_{i,j}}\s{n}{i,j}\s{n}{i,j}^*. 
\]
The unitary $U_{i,j}\in C\big([0,1],M_2(\tcT_{A,B})\big)$ 
defined by 
\[
U_{i,j}(t):=\mat{\sqrt{1-t^2}}
{-t}{t}{\sqrt{1-t^2}}
\widetilde{W_{i,j}(t)}
\mat{\sqrt{1-t^2}}{t}{-t}{\sqrt{1-t^2}}
\mat{\tu_{i,j}^{-1}}{0}{0}{1}
\]
for $t\in [0,1]$ 
satisfies the desired conditions. 
\end{proof}

Note that although $\{U_{i,j}(1)\}_{(i,j)\in \Omega_A}$ 
commute with each others, 
$\{U_{i,j}\}_{(i,j)\in \Omega_A}$ need not. 
However we can show that these unitaries 
commute with each others ``up to homotopy'' 
without changing the values at the two ends. 

\begin{lemma}\label{Lem:a_s}
Let $\cA$ be a unital \CA , 
and $u,v$ be unitaries in $\cA$. 
Suppose that $u,v$ can be written as 
$u=\prod_{i=1}^K\tw_i^{k_i}$, $v=\prod_{i=1}^K\tw_i^{l_i}$ 
for $k_i,l_i\in\Z$ and 
a family $\{w_i\}_{i=1}^K$ of mutually orthogonal partial unitaries. 
Then there exists a homotopy 
$[-1,1]\ni s\mapsto a_s\in \cA$ with $a_{-1}=v-1$, $a_1=u(v-1)$ 
such that $a_s$ commutes with $u$ and $v$, $u^*a_s+va_{-s}^*=0$ 
and $a_sa_s^*=a_s^*a_s=2-v-v^*$ hold for all $s\in [-1,1]$. 
\end{lemma}

\begin{proof}
We may assume that $u=w^k$ and $v=w^l$ 
for a unitary $w\in \cA$ and $k,l\in\Z$. 
We set $a_s=f_s(w)$ for $s\in [-1,1]$ 
where $f_s\in C(\T)$ is defined by 
\[
f_s\big(e^{2\pi \sqrt{-1}t}\big)
:=e^{(s+1)\pi \sqrt{-1}kt}\big(e^{2\pi \sqrt{-1}lt}-1\big)
\]
for $t\in [0,1]$. 
Now it is routine to check that 
the continuous path $[-1,1]\ni s\mapsto a_s\in \cA$ 
satisfies the desired conditions. 
\end{proof}

\begin{lemma}\label{Lem:homotopy}
Let $\cA$ be a unital \CA . 
Let $u,v\in C([0,1],\cA)$ be 
two continuous paths of unitaries 
such that $u(0)=v(0)=1$ and 
$u(1),v(1)\in\cA$ satisfy the assumption in Lemma~\ref{Lem:a_s}. 
Then there exists a homotopy $[-1,1]\ni s\mapsto W_s\in C([0,1],M_2(\cA))$ 
of unitaries such that 
\[
W_{-1}(t)=\mat{u(t)v(t)}{0}{0}{1},\qquad 
W_1(t)=\mat{v(t)u(t)}{0}{0}{1}
\]
for all $t\in [0,1]$ and 
\[
W_s(0)=\mat{1}{0}{0}{1},\qquad 
W_s(1)=\mat{u(1)v(1)}{0}{0}{1}=\mat{v(1)u(1)}{0}{0}{1}
\]
for all $s\in [-1,1]$. 
\end{lemma}

\begin{proof}
For $s\in [0,1]$, 
we define unitaries $W'_{-s},W'_s\in C([0,1],M_2(\cA))$ by 
\begin{align*}
W'_{-s}(t)&=
\mat{u(t)}{0}{0}{1}
\mat{s}{-\sqrt{1-s^2}}{\sqrt{1-s^2}}{s}
\mat{v(t)}{0}{0}{1}
\mat{s}{\sqrt{1-s^2}}{-\sqrt{1-s^2}}{s},\\ 
W'_s(t)&=
\mat{s}{-\sqrt{1-s^2}}{\sqrt{1-s^2}}{s}
\mat{v(t)}{0}{0}{1}
\mat{s}{\sqrt{1-s^2}}{-\sqrt{1-s^2}}{s}
\mat{u(t)}{0}{0}{1}. 
\end{align*}
Then the two definitions of $W'_0$ coincide, 
and we get a homotopy $[-1,1]\ni s\mapsto W'_s\in C([0,1],M_2(\cA))$ 
of unitaries. 
This homotopy satisfies the desired conditions 
except the last condition on $t=1$. 
Since $u(1),v(1)\in\cA$ satisfy the assumption in Lemma~\ref{Lem:a_s}, 
we can find a homotopy 
$[-1,1]\ni s\mapsto a_s\in \cA$ with $a_{-1}=v(1)-1$, $a_1=u(1)(v(1)-1)$ 
such that $a_s$ commutes with $u(1)$ and $v(1)$, 
$u(1)^*a_s+v(1)a_{-s}^*=0$ 
and $a_sa_s^*=a_s^*a_s=2-v(1)-v(1)^*$ hold for all $s\in [-1,1]$. 
We define a homotopy 
$[-1,1]\ni s\mapsto W''_s\in C([0,1],M_2(\cA))$ of unitaries by 
\[
W''_s(t)=
\mat{u(1)\big(t^2v(1)+(1-t^2)\big)}
{t\sqrt{1-t^2}a_{-s}}
{t\sqrt{1-t^2}a_s}
{t^2+(1-t^2)v(1)}
\]
for $s\in [-1,1]$ and $t\in [0,1]$. 
Then this homotopy satisfies 
\[
W''_{s}(0)=\mat{u(1)}{0}{0}{v(1)}=W'_{0}(1),\qquad 
W''_s(1)=\mat{u(1)v(1)}{0}{0}{1}
\] 
for all $s\in [-1,1]$, 
and $W''_{-1}(t)=W'_{-t}(1)$ 
and $W''_{1}(t)=W'_{t}(1)$ for $t\in [0,1]$. 
We define a homotopy $[-1,1]\ni s\mapsto W_s\in C([0,1],M_2(\cA))$ 
of unitaries by 
\begin{center}
\begin{minipage}{10.5cm}
$W_s(t)=
\begin{cases}
W'_s\big(\big(2/(1+|s|)\big)t\big)
& \text{for $0\leq t\leq (1+|s|)/2$,}\\
W''_{s/(2t-1)}(2t-1)
& \text{for $(1+|s|)/2< t\leq 1$. }
\end{cases}$
\end{minipage}
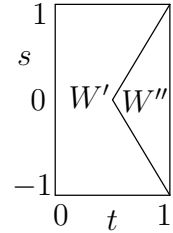
\begin{minipage}{2.5cm}
$\unitlength 0.1in
\begin{picture}(  9.0000, 10.8000)(  3.0000,-12.4500)
\special{pn 8}%
\special{pa 600 200}%
\special{pa 1200 200}%
\special{fp}%
\special{pa 1200 200}%
\special{pa 1200 1200}%
\special{fp}%
\special{pa 1200 1200}%
\special{pa 600 1200}%
\special{fp}%
\special{pa 600 1200}%
\special{pa 600 200}%
\special{fp}%
\special{pa 1200 200}%
\special{pa 900 700}%
\special{fp}%
\special{pa 900 700}%
\special{pa 1200 1200}%
\special{fp}%
\put(7.8000,-6.8000){\makebox(0,0){$W'$}}%
\put(10.7000,-7.0000){\makebox(0,0){$W''$}}%
\put(9.0000,-13.3000){\makebox(0,0){$t$}}%
\put(6.3000,-13.0000){\makebox(0,0){$0$}}%
\put(11.7000,-13.0000){\makebox(0,0){$1$}}%
\put(4.8000,-11.5000){\makebox(0,0){$-1$}}%
\put(5.2000,-2.5000){\makebox(0,0){$1$}}%
\put(5.1000,-7.0000){\makebox(0,0){$0$}}%
\put(4.4000,-4.8000){\makebox(0,0){$s$}}%
\end{picture}$
\end{minipage}
\end{center}
This homotopy satisfies the desired conditions. 
\end{proof}

\begin{remark}
In the lemma above, 
one can weaken the assumptions on $u(1),v(1)\in \cA$. 
However, the assumption that $u(1)v(1)=v(1)u(1)$ 
is too weak to get the conclusion, 
and we can find a counterexample in $\cA=M_2(C(\T^2))$ for example. 
Note that 
the condition that two unitaries $u$ and $v$ commute with each others 
is necessary, but not sufficient 
for the conclusion of Lemma~\ref{Lem:a_s}. 
\end{remark}

Now we construct a $\Gamma_{A,B}$-equivariant splitting map 
$\ker(I-B)\to K_0(\cO_{A,B})$ for 
the surjection $K_0(\cO_{A,B})\to \ker(I-B)$ 
in Proposition~\ref{Prop:exactseqs}. 
The same formulae as the ones in $\cT_{A,B}$ 
define $u_{i,j}\in\cO_{A,B}$ and 
$U_{i,j}\in C\big([0,1],M_2(\tcO_{A,B})\big)$ 
which are nothing but the images 
under the surjection $\cT_{A,B}\to \cO_{A,B}$. 
For $f=\sum_{i=1}^N n_ie_i\in \ker(I-B)$, 
we define a unitary 
\[
U_f:=\prod_{(i,j)\in\Omega_A} U_{i,j}^{n_iB_{i,j}}
\in C\big([0,1],M_2(\tcO_{A,B})\big). 
\]
Then by the computation in Lemma~\ref{Lem:uu^n} 
we have $U_f(0)=U_f(1)=1$. 
Thus this unitary defines an element 
$[U_f]\in K_0(\cO_{A,B})$. 
Although the definition of 
a unitary $U_f\in C\big([0,1],M_2(\tcO_{A,B})\big)$ 
depends on the order of the product above, 
Lemma~\ref{Lem:homotopy} shows that 
the element $[U_f]\in K_0(\cO_{A,B})$ does not depend 
because $U_{i,j}(0)=1$ for all $(i,j)\in\Omega_A$ and 
$\{U_{i,j}(1)\}_{(i,j)\in\Omega_A}$ satisfies 
the assumption in Lemma~\ref{Lem:a_s}. 
From this fact, 
we see that the map 
\[
\ker(I-B)\ni f \mapsto [U_f]\in K_0(\cO_{A,B})
\]
is a well-defined group homomorphism. 
We can also see that this homomorphism is $\Gamma_{A,B}$-equivariant 
by Lemma~\ref{Lem:DefU}. 
We will see that this homomorphism is a splitting map 
for the surjection $K_0(\cO_{A,B})\to \ker(I-B)$. 
The same formula as $U_f$ defines a lifting of $U_f$ 
in $C\big([0,1],M_2(\tcT_{A,B})\big)$ whose value at $1$ is 
$\prod_{i=1}^N (\tu'_i)^{n_i}$ by Lemma~\ref{Lem:uu^n}. 
We see that $[\prod_{i=1}^N (\tu'_i)^{n_i}]\in K_1(\cJ_{A,B})$ 
corresponds to $f=\sum_{i=1}^N n_ie_i\in \ker(I-B)\subset \Z^n$ 
via the natural isomorphism $K_1(\cJ_{A,B})\cong \Z^n$. 
Hence the map $f \mapsto [U_f]$ is a splitting map 
for the surjection $K_0(\cO_{A,B})\to \ker(I-B)$. 

In a similar way as above, 
one can show that the other surjection 
$K_1(\cO_{A,B})\to \ker(I-A)$ in Proposition~\ref{Prop:exactseqs} 
also has a $\Gamma_{A,B}$-equivariant splitting map. 
However for this surjection, 
we can use the well-known computations of $K$-groups of 
Cuntz-Krieger algebras $\cO_A$ as follows (cf.\ \cite{EL2}). 

\begin{proposition}
Let $A,B\in M_N(\Z)$ satisfy the condition (0). 
Then $K_0(\cO_A)$ and $K_1(\cO_A)$ 
are isomorphic to $\coker(I-A)$ and $\ker(I-A)$
as $\Gamma_{A,B}$-modules respectively, 
and the inclusion $\cO_A\hookrightarrow \cO_{A,B}$ 
induces the injection $\coker(I-A)\to K_0(\cO_{A,B})$ 
in Proposition~\ref{Prop:exactseqs} 
and a $\Gamma_{A,B}$-equivariant map 
$\ker(I-A)\to K_1(\cO_{A,B})$ 
which is a splitting map for the surjection 
$K_1(\cO_{A,B})\to \ker(I-A)$ in Proposition~\ref{Prop:exactseqs}. 
\end{proposition}

\begin{proof}
The proof is very similar to the one of Proposition~\ref{Prop:exactseqs}. 

Let $\cT_A\subset \cT_{A,B}$ be the \Csa 
generated by $\{p_i\}_{i=1}^N$ 
and $\{\s{n}{i,j}\}_{(i,j)\in \Omega_{A},n\in\{1,2,\ldots,A_{i,j}\}}$. 
We set $\cJ_A:=\cT_A\cap \cJ_{A,B}$. 
Then we get the following $\Gamma_{A,B}$-equivariant 
commutative diagram with exact rows 
\[
\begin{CD}
0 @>>> \cJ_{A} @>>> \cT_{A} @>>> \cO_{A} @>>> 0\phantom{.}\\
@. @VVV @VVV @VVV @. \\
0 @>>> \cJ_{A,B} @>>> \cT_{A,B} @>>> \cO_{A,B} @>>> 0.
\end{CD}
\]
Let $\cB_N\subset \cA_N$ be the \Csa 
generated by $\{p_i\}_{i=1}^N$. 
Thus $\cB_N\cong \C^N$. 
In a similar way as Lemma~\ref{Lem:cA} and Lemma~\ref{Lem:cA'}, 
we can show that the restrictions $\cB_N\to \cT_{A}$ 
and $\cB_N\to \cJ_{A}$ of \shoms in the two lemmas 
are $KK$-equivalences. 
From these facts and an easy diagram chasing, 
we get the conclusions. 
\end{proof}

This finishes the proofs of the existence 
of the $\Gamma_{A,B}$-equivariant splitting maps 
for the two surjections in Proposition~\ref{Prop:exactseqs}, 
and hence finishes the proof of Proposition~\ref{Prop:exactseqs}.

\end{document}